%% file: maximal_defect_dihedral.tex
\documentclass[11pt,a4paper,twoside]{amsart}

\usepackage{amssymb,amsthm}
\usepackage{wasysym}
\usepackage{array}
\usepackage{enumitem}
\usepackage{tikz}                       
\usetikzlibrary{angles}
\usetikzlibrary{decorations.markings,calc}
\usetikzlibrary{cd}                   
\usepackage{longtable}
\usepackage{mathrsfs}
\usepackage{hhline}
\usepackage{blkarray}        
\usepackage{lscape}
\usepackage{rotating}
\usepackage{arydshln}

\usepackage{comment}
\usepackage{nicematrix}
\usepackage{doi}
\usepackage{titlesec}
\usepackage{orcidlink}

\usepackage{hyperref}                
\hypersetup{                     
pdffitwindow=true,
colorlinks=true,
citecolor=black,
filecolor=black,
linkcolor=black,
urlcolor=black,
hypertexnames=true
}
\usepackage[capitalise,nameinlink]{cleveref}            

\usepackage[textwidth=22mm]{todonotes}
\reversemarginpar

\input xy
\xyoption{all}

\entrymodifiers={!!<0pt,0.7ex>+}  

\setlist{itemsep=1mm, topsep=2mm, leftmargin=10mm}

\newcommand{\IZ}{{\mathbb{Z}}}

\newcommand{\fA}{{\mathfrak{A}}} 
\newcommand{\fS}{{\mathfrak{S}}}  

\newcommand{\cO}{{\mathcal{O}}}

\newcommand{\bA}{{\mathbf{A}}}
\newcommand{\bB}{{\mathbf{B}}}
\newcommand{\bC}{{\mathbf{C}}}
\newcommand{\bD}{{\mathbf{D}}}
  
\newcommand{\bb}{{\mathbf{b}}}
\newcommand{\bc}{{\mathbf{c}}}
\newcommand{\calF}{\mathcal{F}}

\DeclareMathOperator{\Id}{Id}               
\DeclareMathOperator{\End}{End}               
\DeclareMathOperator{\Hom}{Hom}             
\DeclareMathOperator{\Aut}{Aut}                  
\DeclareMathOperator{\Ext}{Ext}                   %
\DeclareMathOperator{\Syl}{Syl}                  
\DeclareMathOperator{\Bl}{Bl}                   
\DeclareMathOperator{\Res}{Res}               
\DeclareMathOperator{\Ind}{Ind}                  
\DeclareMathOperator{\Inf}{Inf}                    

\DeclareMathOperator{\Irr}{Irr}

\DeclareMathOperator{\Sc}{Sc}
\DeclareMathOperator{\Br}{Br}

\DeclareMathOperator{\GL}{GL}

\DeclareMathOperator{\PSL}{PSL}
\DeclareMathOperator{\PGL}{PGL}
\DeclareMathOperator{\PGammaL}{P{\Gamma}L}

\newcommand{\lconj}[2]{\,^{#1}\!#2}                     

\let\wh=\widehat

\newtheorem{thm}{Theorem}[section]
\newtheorem{lem}[thm]{Lemma}

\newtheorem{cor}[thm]{Corollary}
\newtheorem{prop}[thm]{Proposition}

\newtheoremstyle{defnew}{4ex}{}{}{}{\bf}{.}{.5em}{}
\theoremstyle{defnew}

\newtheorem{nota}[thm]{Notation}
\newtheorem{rem}[thm]{Remark}

\theoremstyle{remark}

\theoremstyle{plain}

\titleformat{\section}
  {\normalfont\scshape\centering}
  {\thesection.}
  {0.2em}
  {}

\titlespacing*{\section}{0pt}{6ex}{2ex}

\titleformat{\subsection}[runin]
  {\normalfont\bfseries}
  {\thesubsection }
  {0.5em}
  {}
  [.\enspace]

\titlespacing*{\subsection}
  {0pt}
  {0.5\baselineskip plus 0.7\baselineskip}
  {0pt}

\begin{document}


\title[Maximal defect $2$-blocks with dihedral defect groups up to splendid Morita equivalence]{Maximal defect $2$-blocks with dihedral defect groups up to splendid Morita equivalence}
	
\date{\today}

\author{{Shigeo Koshitani\,\raisebox{0.6ex}{\orcidlink{0000-0002-0703-0363}}, Caroline  Lassueur\,\raisebox{0.6ex}{\orcidlink{0000-0001-9031-0587}}, {\.I}pek Tuvay\,\raisebox{0.6ex}{\orcidlink{0000-0002-7427-9310}}, Deniz Y{\i}lmaz}\,\raisebox{0.6ex}{\orcidlink{0000-0002-4245-8154}}
}
\address{{\sc Shigeo Koshitani}, Department of Mathematics and Informatics,
Graduate School of Science, Chiba
University, 1-33 Yayoi-cho, Inage-ku, Chiba 263-8522, Japan}
\email{koshitan@math.s.chiba-u.ac.jp} 
\address{{\sc Caroline Lassueur}, Leibniz Universit\"at Hannover,
Institut f\"ur Algebra, Zahlentheorie und Diskrete Mathematik,
Welfengarten 1, 30167 Hannover, Germany}
\email{lassueur@mathematik.uni-kl.de}
\address{ {\sc {\.I}pek Tuvay},  Mimar Sinan Fine Arts University, Department of Mathematics, 34380 Bomonti, İstanbul, Turkey}
\email{ipek.tuvay@msgsu.edu.tr}
\address{{\sc Deniz Y{\i}lmaz},  Department of Mathematics, 
Bilkent University
Ankara, Turkey}
\email{d.yilmaz@bilkent.edu.tr}



\makeatletter
\@namedef{subjclassname@2020}{%
		\textup{2020} Mathematics Subject Classification}
	\makeatother

	\begin{abstract}
        We classify, up to splendid Morita equivalence, all maximal defect blocks whose defect groups are dihedral $2$-groups of order at least $8$. The classification holds for blocks over an algebraically closed field as well as for those over a suitable discrete valuation ring of characteristic zero. As a consequence, Puig’s conjecture holds for this infinite family of blocks.\par
        Beyond the classification result itself, a central feature of the present work is the variety of methods required in its proof. Our arguments combine classical techniques of block theory with a broad range of tools and techniques that have emerged over the past twenty-five years.
         Amongst others these include in a crucial manner the Gorenstein-Walter  classification of finite groups with dihedral Sylow $2$-subgroups, Linckelmann's glueing of local Morita equivalences to obtain a stable equivalence of Morita type, the modern versions of Fong's 1st and 2nd reduction theorems in terms of splendid and endo-permutation source Morita equivalences, the structure of the Dade group of endo-permutation modules, Boltje-Perepelitsky's extended theory of Brauer pairs for $p$-permutation modules, the fusion systems of the blocks involved, or the concrete knowledge of their decomposition matrices.
	\end{abstract}

\thanks{
The authors gratefully acknowledge the generous support and hospitality of the Mathematisches Forschungsinstitut Oberwolfach, where this project  developed in January 2026 as part of the Oberwolfach Research Fellowship Programme for research in small groups.
}
	
\maketitle

\pagestyle{myheadings}
\markboth{S. Koshitani, C. Lassueur, {\.I}. Tuvay, D. Y{\i}lmaz}{Maximal defect $2$-blocks with dihedral defect groups}


\section{Introduction}\label{sec:intro}%

Let $p$ be a prime number, let $k$ be an algebraically closed field of characteristic~$p$ and let $G$ be a finite group. An indecomposable direct algebra factor of the group algebra~$kG$ is called a \emph{block} (or \emph{block algebra}, or a \emph{$p$-block} of  $G$). The present work is part of the general programme of classifying blocks up to certain categorical equivalences, which preserve various representation-theoretic invariants.
\par
Two notions of equivalence play a central role in this programme:
the standard notion of \emph{Morita equivalence} of finite-dimensional $k$-algebras and the stronger notion of \emph{splendid Morita equivalence}. 
Donovan's Conjecture and Puig's Finiteness Conjecture predict, respectively, that there are only finitely many Morita equivalence classes and splendid Morita equivalence classes $p$-blocks with a given isomorphism type of defect groups.
A proof of Donovan's Conjecture has recently been claimed by OpenAI~\cite{OpenAIDonovan, OpenAIIntegralDonovan}.
Beyond these finiteness statements, a natural, and much more involved and deeper task than their numerical counting is to classify the equivalence classes explicitly.
\par
The fast-growing database of blocks maintained by Charles Eaton \cite{DATABASE_CHARLES_EATON} provides an excellent overview of the progress made towards classifying blocks up to Morita equivalence.
Much less is known about their classification up to splendid Morita equivalence. Puig’s Finiteness Conjecture  holds for blocks with cyclic defect groups, for any prime number $p$ \cite{LinckCyclic96}, and for blocks with Klein-four defect groups when $p=2$ \cite{CELK11}. 
In the latter case, \cite{CELK11} also provides us with a complete classification: there are three splendid Morita equivalence classes. For cyclic defect groups, however, a complete classification is not yet known, even up to Morita equivalence. 
A good summary of the known/unknown embedded Brauer trees can be found in \cite[Chapter 5]{CravenRepBook}, while   \cite{HissLassueur24, HissLassueur25, HissLassueur26} investigate one of the further invariants needed to determine the source algebra of such blocks.
\par
These considerations illustrate the difficulty of classifying blocks with a given defect group up to splendid Morita equivalence, even when their representation theory is well-understood.
Nevertheless, further finiteness results in the direction of Puig's Finiteness Conjecture have been obtained by restricting the class
of finite groups, respectively the class of blocks, considered. These include \cite{PuigNilpotentBlocks88} for all blocks of $p$-soluble
groups, \cite{Pui94, Kes96} for all blocks of the symmetric groups and their double covers, \cite{Kes96, Kes02}
for all blocks of the alternating groups and their double covers, \cite{Kes00} for all blocks of Weyl groups,
\cite{HK05, HK00, Kes01} for the unipotent blocks of the finite classical groups, and \cite{KoLa19,KL20,KLS21} for principal blocks of tame representation type. 
\par
In this spirit, a natural next step is to investigate arbitrary $2$-blocks of tame representation type, i.e. with dihedral, semi-dihedral or generalized quaternion defect groups. In this article, we take a first step in this direction by considering tame blocks of maximal defect with defect groups isomorphic to a dihedral $2$-group $D^{}_{2^n}$ of order $2^n$ with $n\geq 3$.  Our main result shows that every such block is splendidly Morita equivalent to a principal block.


\begin{thm}\label{thm:mainthm_intro}
   Let $k$ be an algebraically closed field of characteristic $2$ and let $G$ be a finite group with a dihedral Sylow $2$-subgroup of order $2^n$ $(n\geq 3)$. Then, any block $\bB$ of $kG$ with maximal defect is splendidly Morita equivalent to precisely one of the following blocks: 
    \begin{enumerate}[label=\rm(\alph*)]
       \item $\bB_0(kD_{2^n})$; 
       \item $\bB_0(k\fA_7)$  (provided $n=3$); 
       \item $\bB_0(k\!\PSL_2(q))$ where $q$ is an odd prime power such that  $(q-1)_2=2^{n}$;
       \item $\bB_0(k\!\PSL_2(q))$ where $q$ is an odd prime power such that  $(q+1)_2=2^{n}$;
       \item $\bB_0(k\!\PGL_2(q))$ where $q$ is an odd prime power such that  $(q-1)_2=2^{n-1}$;
       \item $\bB_0(k\!\PGL_2(q))$ where $q$ is an odd prime power such that  $(q+1)_2=2^{n-1}$.
    \end{enumerate}
     In cases {\rm(c)--(f)}, any choice of q satisfying the indicated condition yields a representative of the corresponding splendid Morita equivalence class.   
 \end{thm}  

We recall that Erdmann classified tame $k$-algebras up to Morita equivalence through an explicit list of families of basic algebras, defined by quivers and relations; see \cite{ErdmannTameBook}. This classification was recently refined by Macgregor \cite{MacGregor22} using the classification of the finite
simple groups, giving a complete list of Morita equivalence classes of blocks of group algebras with dihedral defect groups.
Theorem~\ref{thm:mainthm_intro} shows that the list obtained in
\cite{MacGregor22} does not split further up to splendid Morita equivalence if we restrict ourselves to blocks of maximal defect.

Now, if we consider a complete discrete valuation ring $\cO$  of characteristic zero with residue field $\cO/J(\cO)=k$ (our previously fixed algebraically closed field of characteristic~$2$), then the results of Eisele in \cite{Eiselethesis} prove that Erdmann's and Macgregor's classification can be lifted to $\cO$.  Adding the splendid condition, 
Theorem~\ref{thm:mainthm_intro} can be lifted to $\cO$ as a splendid Morita equivalence given by tensoring with a $2$-permutation bimodule, which is liftable.
Moreover, Theorem~\ref{thm:mainthm_intro}   directly implies that Puig's finiteness conjecture holds when we restrict ourselves to the class of all blocks of all finite groups with defect groups isomorphic to a fixed dihedral $2$-group and maximal defect.

\begin{cor}\label{cor:intro_PuigsConj}%
    For any $n\geq 3$, there are only finitely many splendid Morita equivalence classes of blocks of the group algebra of finite groups, over $k$ and over $\cO$, with defect groups isomorphic to $D_{2^n}$ and maximal defect. 
\end{cor}

We emphasize that beyond the classification result of Theorem~\ref{thm:mainthm_intro} itself, a significant aspect of the present work lies in the wide range of methods and techniques that are needed for its proof. The arguments combine many classical tools from group and block theory with many more recent results and methods that have emerged over the past 25 years. See the abstract.
Moreover, the dihedral case considered here will also serve as a prototype for the treatment of the remaining tame cases, that is, the cases of generalized quaternion and semi-dihedral defect group, which we will treat in forthcoming papers. The recent work of Eaton, Eisele, Kessar, Linckelmann and Schaeffer Fry \cite{EEKLSF26} classifying $2$-blocks with generalized quaternion defect groups up to Morita equivalence will be of substantial help in the treatment of this case.

We also emphasize that the main difficulty at the heart of the proof lies in the problem of establishing splendid Morita equivalences between the faithful \(2\)-blocks of maximal defect of \(3.\fA_6\) and \(3.\fA_7\) and suitable principal blocks.

The paper is organized as follows. To begin with, in Section~\ref{sec:prelim} and~\ref{sec:blocksWithDihedralDefectGroups},  we introduce the background material which is essential for the proof of Theorem~\ref{thm:mainthm_intro}. In Section~\ref{sec:4}, we perform a reduction to $G/O_{2^\prime}(G)$ and its central extensions by a certain odd cyclic group. In Section~\ref{sec:3A7} and~\ref{sec:3A6}, we deal with four isolated cases arising from the exceptional covering groups $3.\fA_6$ and $3.\fA_7$.
Finally, we prove our main results in Section~\ref{sec:proofs}. 
Furthermore, in  Appendix~\ref{app:A} we give a proof that Fong's second reduction theorem gives a Morita equivalence with endo-permutation source and diagonal vertex.



\input{section_2_notation_and_quoted_results}

\input{section_3_blocks_with_dihedral_defect_grps_prelim_res}

\input{section_4_reduction_to_GmodO2}

\input{sections_5_6_and_7}



\bibliographystyle{aomalpha}
\bibliography{biblio.bib}
\bigskip



\input{appendix}


\end{document}

%% file: section_2_notation_and_quoted_results.tex
\section{Notation, prerequisites and quoted results}\label{sec:prelim}%

\subsection{General notation}
Throughout this manuscript, unless otherwise specified, we adopt the following notation and conventions.
\par
All groups considered are finite, and $G$ and $H$ always denote finite groups. We let $k$ be an algebraically closed field of characteristic~$p>0$.  From Section~\ref{sec:blocksWithDihedralDefectGroups} on, we assume that $p=2$.
When character-theoretic  arguments are needed, we let  $\cO$ be a complete discrete valuation ring with residue field $\cO/J(\cO)=k$ and a quotient field $K$ of characteristic zero, which is
assumed to be large enough for the finite groups under consideration and their subgroups. 
\par 
For $R\in\{k,\cO\}$ all modules over a given $R$-algebra $\bA$ that is free of finite rank over \(R\), are assumed to be finitely generated left
modules and we denote by $\bA-\mathsf{mod}$ the corresponding module category.  
The set of all blocks of $\bA$ (i.e. its indecomposable algebra factors) is denoted $\Bl(\bA)$.
We write $\bB_0(RG)$ for the principal block of $RG$. 
Given $\bB\in\Bl(RG)$, we write  $e_{\bB}$  for the block idempotent of $RG$ corresponding to $\bB$, and 
if $D\leq G$ is a defect group of $\bB$, then the integer $d$ such that $|D|=p^{d}$ is called the defect of $\bB$. We say that $\bB$ is a block of \emph{maximal defect} provided its defect groups are the Sylow $p$-subgroups of $G$, i.e.  $|G|_p=p^d$, and 
we let $\Bl^{\mathrm{max}}(RG)$ denote  the subset of $\Bl(RG)$ consisting of all blocks of maximal defect of $RG$. If $N\unlhd G$ and $\bc\in\Bl(RN)$, then  $\Bl(RG\mid \bc):=\{\bB\in\Bl(RG)\mid \bB\text{ covers }\bc\}$.

\par
The canonical algebra homomorphism $\cO G\longrightarrow kG$ induces a defect group preserving bijection between $\Bl(\cO G)$ and $\Bl(kG)$, mapping a block $\bB\in \Bl(\cO G)$ to $k\otimes_{\cO}\bB$. In this context, we write $\Irr(\bB)$,  and also simply $\Irr(k\otimes_{\cO}\bB)$ when working over $k$, 
for the sets of all
irreducible 
$K$-characters 
of $G$ belonging to $\bB$. We write $\Irr(G)$ for the set of all irreducible $K$-characters of~$G$.
\par
In contrast, if $\bB\in\Bl(kG)$, then $\Irr_k(\bB)$ denotes the set of all simple $\bB$-modules up to isomorphism. 
We denote by $k_G$, or simply $k$,  the trivial $kG$-module. 
For $H\leq G$ we denote by $\Res^G_H$ (resp.~$\Ind_H^G$) the restriction (resp. the induction) of $kG$-modules and characters. Similarly, for $N\unlhd G$, we denote by $\Inf_{G/N}^{G}$ the inflation of $kG$-modules and characters. Given a $p$-subgroup $Q\leq G$ and a $kG$-module $M$  we denote by $M[Q]$ the Brauer construction of $M$ with respect to $Q$, and by $\Br^G_Q$ (or simply $\Br_Q$ when clear from the context)  the Brauer homomorphism  with respect to $Q$, i.e. the surjective algebra homomorphism $\Br^G_Q: (kG)^Q\longrightarrow kC_G(Q)$
mapping $C_G(Q)$ to itself by the identity map and mapping all the other group
elements to zero. The projective cover of a simple $kG$-module $S$ is written $P(S)$. If $S_1,S_2$ are simple $kG$-modules, then $U(S_1,S_2)$ denotes a uniserial $kG$-module with head $S_1$ and socle $S_2$. 
\par
Again, given $R\in\{\mathcal{O},k\}$, an $RG$-lattice $M$ is a \emph{permutation} $RG$-lattice if it admits a $G$-invariant $R$-basis, whereas it is a $p$-\emph{permutation} $RG$-lattice if it is a direct summand of a permutation $RG$-lattice. Indecomposable $p$-permutation $RG$-lattices are also called \emph{trivial source} $RG$-lattices. 
It is well known that any $p$-permutation $kG$-module $M$ is liftable to an $\cO G$-lattice. In general, such modules afford several lifts, but there is a unique one amongst these which is a $p$-permutation $\cO G$-lattice. We denote this lift by $\wh{M}$ and we let $\chi_{\wh{M}}$ be the character afforded by $K\otimes_{\cO}\wh{M}$.
For any $p$-subgroup $Q$ of $G$, the Scott module with vertex $Q$, written $\Sc(G,Q)$, is always a trivial source module. 
\par
Given an integer $n\geq 2$, 
$C_n$ denotes the cyclic group of order $n$, $\fA_n$ and $\fS_n$ the alternating and the symmetric group on $n$ letters, respectively, and $D_{2^n}$ the dihedral group of order $2^n$. 
Furthermore,   ${\mathrm{Syl}}_p(G)$ denotes  the set of all Sylow $p$-subgroups of $G$, and for $P\in\Syl_p(G)$ we let $\mathcal F_P(G)$ be the fusion system of $G$ on $P$. 
If $\varphi:Q\longrightarrow P$ is an isomorphism between subgroups $Q\leq H$ and $P\leq G$, then $\Delta(P,\varphi,Q):=\{(\varphi(y),y)\in P\times Q\mid y\in Q\}$
is a twisted diagonal subgroup of $G\times H$, and if $G=H$ and $P=Q$ then $\Delta(P):=\Delta(P,\Id,P)$ is the diagonal embedding of $P$ in $G\times G$. Further, we denote the (external) central product of two finite groups $G_1$ and $G_2$ identified along a group isomorphism $\varphi:Z_1\longrightarrow Z_2$ with $Z_i\leq Z(G_i)$ for $i=1,2$ by $G_1\ast G_2$, so 
$G_1\ast G_2=(G_1\times G_2)/Z$ with $Z=\{(z,\varphi(z)^{-1})\mid z\in Z_1\}\leq Z(G_1\times G_2)$.
\par
For further standard notation and terminology, we refer the reader to \cite{LinckelmannVol1,LinckelmannVol2,NagaoTsushimaBook,ThevenazBook}.

\vspace{2mm}
\subsection{Definitions: Morita equivalences, stable equivalences of Morita type, splendid Morita equivalences and endo-permutation source Morita equivalences}\label{ssec:Splendid_Morita_equivalences}{\leavevmode}%

\noindent Let $\bA$ and $\bB$ be finite-dimensional $R$-algebras with $R\in\{\cO,k\}$. We recall that a pair $(M,N)$ consisting of an $(\bA,\bB)$-bimodule $M$ and of a $(\bB,\bA)$-bimodule $N$ is said to \emph{induce a Morita
equivalence between $\bA$ and $\bB$} if $M$ and $N$ are finitely generated projective as left
and right modules, and if there are isomorphisms $M\otimes_{\bB} N\cong \bA$ and $N\otimes_{\bA}M\cong \bB$ of $(\bA,\bA)$-bimodules and $(\bB,\bB)$-bimodules respectively. In this case, the functors $M\otimes_{\bB}–:\bB-\mathrm{mod}\longrightarrow \bA-\mathsf{mod}$ and $N\otimes_{\bA}–:\bA-\mathrm{mod}\longrightarrow \bB-\mathsf{mod}$ are inverse equivalences of categories, and $\bA$ and $\bB$ are called \emph{Morita equivalent}.  Moreover, as bimodules, $N\cong  \Hom_{\bA}(M,\bA)$ and $M\cong \Hom_{\bB}(N,\bB)$, so we often simply say that $M$ \emph{induces a Morita
equivalence between $\bA$ and $\bB$}.  See~\cite[Theorem~2.8.2]{LinckelmannVol1}. 
Under the weaker assumptions that $M\otimes_{\bB} N\cong \bA\oplus(\text{projective $(\bA,\bA)$-bimodule})$ and $N\otimes_{\bA}M\cong \bB\oplus(\text{projective $(\bB,\bB)$-bimodule})$,  it is said that the pair $(M,N)$ induces a \emph{stable equivalence of Morita type} between $\bA$ and $\bB$. 

\par
Assume now that $G$ and $H$ are finite groups, $\bA\in\Bl(RG)$, and $\bB\in\Bl(RH)$. Then, the $(\bA,\bB)$-bimodule $M$ can be regarded as a left $R[G\times H]$-module (via the left action given by $(g,h)\cdot m:=g\cdot m\cdot h^{-1}$ for any $(g,h)\in G\times H$ and any $m\in M$) and the $(\bB,\bA)$-bimodule $N$ as a left $R[H\times G]$-module. The Morita equivalence is called a \emph{splendid Morita equivalence} (respectively, an \emph{endo-permutation source Morita equivalence}) if $M$ and $N$, seen as left modules, are trivial source modules (respectively, have an endo-permutation source). Clearly, any splendid Morita equivalence is an endo-permutation source  Morita equivalence. If there is an endo-permutation source Morita equivalence between $\bA$ and $\bB$, then they  have isomorphic defect groups. After identifying defect groups of $\bA$ and $\bB$ with a $p$-group $D$, the modules $M$ and $N$ have vertex $\Delta D$. A splendid Morita equivalence is sometimes also called a \emph{Puig equivalence} or a \emph{source-algebra equivalence}, since it is equivalent to the existence of an isomorphism between source algebras of $\bA$ and $\bB$ as interior $D$-algebras. See \cite[Proposition~9.7.1, Remark~9.7.2, Definition~9.7.3,  Theorem~9.7.4 and Theorem 9.11.2]{LinckelmannVol2}.
\vspace{2mm}
\subsection{Blocks of direct and central products}\label{ssec:blocks_of_central_producrs}
Let $G$ and $H$ be finite groups as in the previous subsection. 
Recall that there is a canonical $R$-algebra isomorphism $R[G\times H]\longrightarrow RG\otimes_R RH$, mapping $(g,h)\in G\times H$ to $g\otimes h$. Identifying $R[G\times H]$ and $RG\otimes_R RH$ via this isomorphism, each block $\bB \in\Bl(R[G\times H])$ can be written as $\bB=\bB_G\otimes_R\bB_H\in \Bl(RG\otimes_R RH)$ with $\bB_G\in\Bl(RG)$ and $\bB_H\in\Bl(RH)$. 
\par
More generally, if $G\ast H$ is a central product of $G$ and $H$ identified along a central subgroup of $p'$-order, then any block $\bB \in\Bl(R[G\ast H])$ is isomorphic to $\bB_G\otimes_R\bB_H$ where $\bB_G\in\Bl(RG)$ is the unique block of $RG$ covered by $\bB$ and $\bB_H\in\Bl(RH)$ is the unique block of $RH$ covered by $\bB$. See \cite[Proposition~1.5]{Duvel04} (or also \cite[Lemma~7.5]{Sambale14}). More precisely, the block $\bB_G\otimes_R\bB_H$, viewed as a block of \(R[G\times H]\), maps onto $\bB$ under the canonical map $R[G\times H]\to R[G\ast H]$. By
\cite[Lemma~1.4]{Duvel04}, this map induces an algebra isomorphism
\[
        \bB_G\otimes_R\bB_H\cong \bB.
\]
Since the kernel of the canonical homomorphism $G\times H\to G\ast H$ is a $p'$-group, the quotient map identifies the corresponding defect groups. If $D$ denotes the common defect group, then the above map, being induced by a group homomorphism, is compatible with the corresponding interior $D$-algebra structures. Thus, $\bB_G\otimes_R\bB_H$ and $\bB$ have isomorphic source algebras as interior $D$-algebras. In particular, they are splendidly Morita equivalent.

\begin{lem}\label{lem:SplendidMoritaEquivalence_direct_and_central_product}%
Let $G\ast H$ be a central product of finite groups identified along a central subgroup of $p'$-order (possibly trivial).
Let $\bB \in\Bl(R[G\ast H])$ and let $\bB_G\in\Bl(RG)$ and $\bB_H\in\Bl(RH)$ be the unique blocks covered by $\bB$, as above. If $\bB_G$ is a block of defect zero, then $\bB$ and $\bB_H$ are splendidly Morita equivalent. 
\end{lem}

\begin{proof}
By the above discussions, the blocks $\bB$ and $\bB_G\otimes_R\bB_H$ are splendidly Morita equivalent. By \cite[p.\,341, line -9]{ThevenazBook}, $R$ is a source algebra of $\bB_G$, thus the assertion follows from \cite[Lemma 2.6(i)]{KMN11}.
\end{proof}

\vspace{2mm}
\subsection{Groups with dihedral Sylow \texorpdfstring{$2$}{2}-subgroups}\label{ssec:dihedralSylow2}%
Let $G$ be a finite group with a dihedral Sylow $2$-subgroup $P$ of order $2^n\geq 8$. Then, Gorenstein and Walter~\cite[Theorem~1]{GorensteinWalter} proved 
that $\overline{G}:=G/O_{2'}(G)$ is isomorphic to $P$, the alternating group $\fA_7$,  or a subgroup of $\PGammaL_2(q)$ containing $\PSL_2(q)$, 
where $q:=r^f$ is a power of an odd prime $r$ and $f\geq 1$ is an integer. Then,  
$\PGammaL_2(q)\cong \PGL_2(q)\rtimes 
{\mathrm{Gal}}(\mathbb F_{q}/\mathbb F_{r})$,  
where ${\mathrm{Gal}}(\mathbb F_{q}/\mathbb F_{r})$  
is cyclic of order $f$. Thus, for the purposes of our computations, we split the above list of groups into further subcases as follows:
\begin{enumerate}[label={\bf(D\arabic*)}, leftmargin=12mm]
   \item $\overline{G}\cong D_{2^n}$;
   \item $\overline{G}\cong \fA_7$;
   \item  $\overline{G}\cong \PSL_2(q)\rtimes C_d$,  where $q$ is a power of an odd prime such that $(q-1)_2=2^{n}$ or $(q+1)_2=2^{n}$, $d$ is odd and $d\,|\,f$;
   \item $\overline{G}\cong \PGL_2(q)\rtimes C_d$,  where $q$ is a power of an odd prime such that $(q-1)_2=2^{n-1}$ or $(q+1)_2=2^{n-1}$, $d$ is odd and $d\,|\,f$.
\end{enumerate}
Notice that Case {\bf(D3)} happens provided $q\equiv\pm 1\pmod{8}$, whereas Case {\bf(D4)} happens when  $q\equiv\pm 1\pmod{4}$. In both cases, $d$ is odd. 
Moreover, $\PSL_2(9)\cong \fA_6$ and $\PGL_2(9)\cong \fA_6.2_2$ in the ATLAS notation. \\

\noindent Next, recall that the Schur multiplier of a finite group $G$ is the abelian group ${\mathrm{H}}^2(G,\mathbb{C}^\times)$ and there is an isomorphism 
\[
    {\mathrm{H}}^2(G,k^\times)\cong {\mathrm{H}}^2(G,\mathbb{C}^\times)_{2^\prime} \tag{$\ast$}
\]
(see e.g. \cite[2.1.14 Proposition]{KarpilovskySchurMultBook}). For the groups above, this group is as follows.

\begin{lem}\label{lem:SchurMult_2prime}
Assume $p=2$. Let $G$ be a finite group with a dihedral Sylow $2$-subgroup of order $2^n$ with $n\geq 3$ and assume that $O_{2^\prime}(G)=1$. Then, we have 
\[
{\mathrm{H}}^2(G, k^\times)\cong  \begin{cases}
    C_3  & \text{ if } \,  G
\in \{\mathfrak{A}_7, \PSL_2(9)\}, \\
    1  & \text{otherwise}.
\end{cases}
\]
Moreover, if $G
\in \{\mathfrak{A}_7, \PSL_2(9)\}$ and $Z$ is a finite cyclic group of odd order on which $G$ acts trivially and $3\mid |Z|$, then 
\[
{\mathrm{H}}^2(G, Z)\cong {\mathrm{H}}^2(G, k^\times)\cong C_3\,.
\]
\end{lem}

\begin{proof} 
We go through the cases {\bf (D1)} to {\bf (D4)}.  First, if $G=D_{2^n}$, then ${\mathrm{H}}^2(G, k^\times)\cong 1$ as $G$ is a $2$-group; see e.g. \cite[Chapter 2, Theorem~7.4]{NagaoTsushimaBook}.
\par
If $G\in\{\fA_7,\PSL_2(9)\}$, then we read from the ATLAS \cite[p.10 and p.4 resp.]{ATLAS} that ${\mathrm{H}}^2(G,\mathbb{C}^\times)\cong C_6$, so ${\mathrm{H}}^2(G, k^\times)\cong C_3$ by ($\ast$). 
\par
If {\bf (D3)} holds with $d=1$, $G=\PSL_2(q)$  and $q\neq 9$, then ${\mathrm{H}}^2(G,\mathbb{C}^\times)\cong C_2$  by \cite[Chapter~3, \S3 with Table~5]{ATLAS}, hence  ${\mathrm{H}}^2(G, k^\times)\cong 1$ by ($\ast$).
\par
Next, if $G\cong \PGL_2(3)\cong\mathfrak S_4$, then we have ${\mathrm{H}}^2(G, \mathbb{C}^\times)\cong C_2$
 by \cite[2.12.3 Theorem]{KarpilovskySchurMultBook}, hence ${\mathrm{H}}^2(G, k^\times)\cong 1$ by ($\ast$). If $G=\PGL_2(9)$, then we read from the ATLAS \cite[bottom of p.5 together with page xxii]{ATLAS} and ($\ast$) that ${\mathrm{H}}^2(G, k^\times)\cong 1$. 
\par
We go through the remaining groups using \cite[Theorem 3.1(c)]{Jones74}, which states that if $N\unlhd G$ is such that $G/N$ is cyclic, then 
$|{\mathrm{H}}^2(G,\mathbb{C}^\times)|$ divides $|{\mathrm{H}}^2(N,\mathbb{C}^\times)|\cdot |N/[N,N]|$. In all cases, we take $N:=\PSL_2(q)$, which is perfect. If $G$ is of type {\bf (D3)} with $d\neq 1$, then $G/N\cong C_d$, and certainly $q>9$.  If $G$ is of type {\bf (D4)} and $q\notin\{3,9\}$, then $G/N\cong C_{2d}$. In both cases, we obtain that 
\[
|{\mathrm{H}}^2(G,\mathbb{C}^\times)|\,\big|\,|{\mathrm{H}}^2(N,\mathbb{C}^\times)|\cdot |N/[N,N]|=2\cdot 1=2\,,
\]
thus ${\mathrm{H}}^2(G, k^\times)\cong 1$ by ($\ast$). This proves the first claim. The second claim follows from the Universal Coefficient Theorem, which reads  
\[
{\mathrm{H}}^2(G, Z) \cong \Hom_{\IZ}({\mathrm{H}}_2(G, \mathbb{Z}), Z) \oplus \Ext^1_{\mathbb{Z}}({\mathrm{H}}_1(G, \mathbb{Z}), Z)\,.
\]
(See e.g. \cite[Exercise~6.1.5]{WeibelBook}.)
Since $\fA_7$ and $\PSL_2(9)$ are perfect, the Ext part vanishes, because  ${\mathrm{H}}_1(G, \mathbb{Z})\cong G/[G,G]$. Moreover, $ {\mathrm{H}}_2(G,\mathbb{Z})\cong {\mathrm{H}}^2(G,\mathbb{C}^\times)$ and we obtain ${\mathrm{H}}^2(G, Z)\cong \Hom_{\IZ}(C_3, Z)\cong C_3$. 
\end{proof}

\vspace{2mm}
\subsection{Endo-permutation modules and the Dade group}\label{subsec:endo-perm}%
Let $P$ denote a $p$-group. A $kP$-module $M$ is called \emph{endo-permutation} if $\End_k(M)$ is a permutation $kP$-module. Moreover, such an endo-permutation $kP$-module $M$ is called \emph{capped} provided it possesses an indecomposable direct summand with vertex $P$, denoted $\mathrm{cap}(M)$ and called its \emph{cap}. In this case, $\mathrm{cap}(M)$ is unique up to isomorphism, although it may appear with multiplicity greater than one, so $M=\mathrm{cap}(M)^a\oplus X$ with $a\in\mathbb{Z}_{>0}$ and $X$ is a $kP$-submodule of $M$, the indecomposable direct summands of which have vertices of order strictly smaller than $|P|$.   
This allows us to define an equivalence relation: two capped endo-permutation $kP$-modules are equivalent if and only if their caps are isomorphic. The \emph{Dade group} of $P$, denoted $\bD_k(P)$, is the resulting set of equivalence classes (of capped endo-permutation $kP$-modules) with group law given by the tensor product $\otimes_k$\,. By a result of Puig, $\bD_k(P)$ is a finitely generated abelian group. We denote its torsion part by $\bD_k^{tors}(P)$. For a detailed introduction to endo-permutation modules and their properties, we refer the reader to \cite{ThevenazSurvey}. \\

\noindent We need the next lemma linking sources of simple modules to the structure of the Dade group, the proof of which consists of the arguments given in \cite[paragraph after Proposition~2.1]{CELK11}.

\begin{lem}\label{lem:SourcesSimples_p_soluble_groups}\label{lem:sourcesSimpleDihedral}\label{lem:SourcesSimplesSolubleGroups}
Let $G$ be a $p$-soluble group. Let $T$ be a simple $kG$-module, let $Q$ be a vertex of~$T$ and let $U$ be a $kQ$-source of $T$.  Then, the following assertions hold:
    \begin{enumerate}[label={\rm(\alph*)}]
        \item $U$ is a capped endo-permutation $kQ$-module and its class in the Dade group is a torsion element, i.e. $[U]\in\bD_k^{tors}(Q)$; and 
        \item if $p=2$ and  $Q\cong D_{2^n}$ with $n\geq 3$, then $U\cong k$. 
    \end{enumerate}
\end{lem}    
    
\begin{proof}    
    (a) First, $U$ is an endo-permutation module by \cite{PuigNilpotentBlocks88}. Furthermore, by a result of Feit, the simple module $T$ is algebraic; see \cite{Feit80}. Thus, \cite[Proposition~4.4]{BoltjeKuls06} proves that the class of $U$ in the Dade group $\bD_k(Q)$ must be a torsion element.\\
    (b) As $Q$ is a dihedral $2$-group, we have  $\bD_k^{tors}(Q)=\{[k]\}$ (see e.g. \cite[Proposition~6.4(d)]{ThevenazSurvey}), proving that $U$ is the trivial $kQ$-module. 
\end{proof}

\vspace{2mm}
\subsection{Clifford-theoretic results}

We will use the standard reduction techniques due to 
Fong and Reynolds. We review them here in terms of splendid and endo-permutation source Morita equivalences, which is less standard in the literature. These techniques work irrespective of the characteristic, so in this subsection $p$ is an arbitrary prime.

\begin{thm}[Fong{--}Reynolds, 
see {\cite[1.5. Theorem]{KoshitaniKunugiWaki04}}]\label{thm:FongReynolds}%
Let $G$ be a finite group, let $N\unlhd G$, let $\bb\in\Bl(kN)$
and  let $T:= T_G(\bb)$ be the stabilizer of $\bb$ in $G$. Then, block induction induces a bijection
$$\Bl(kT\mid\bb)\longrightarrow \Bl(kG\mid\bb), \bA\mapsto \bA^G\,.$$
Moreover, if $\bB\in\Bl(kG\mid\bb)$ and $\widetilde{\bB}\in\Bl(kT\mid\bb)$ is the preimage of $\bB$ under the latter bijection, then the pair of bimodules $(e^{}_{\bB}\cdot kG\cdot e^{}_{\widetilde{\bB}}, e^{}_{\widetilde{\bB}}\cdot kG\cdot e^{}_{\bB})$ realizes a splendid Morita equivalence between $\bB$ and $\widetilde{\bB}$. 
\end{thm}

\begin{thm}[Fong’s Second Reduction Theorem, see {\cite[Theorem~6.8.13]{LinckelmannVol2}} and  {\cite[ Proposition~3.3]{SalminenSOURCES08}}]\label{thm:Fong2ndReduction}%
   Let $G$ be a finite group and let $N\unlhd G$ be a normal subgroup. Let $\bB\in\Bl(kG)$ be a block  covering a $G$-stable block $\bc\in\Bl(kN)$ and set $S:=kNe_{\bc}$. 
   If $D$ is a defect group of $\bB$ and $D\cap N=1$, then $\bc$ has defect zero and there is an isomorphism of interior $D$-algebras
   \[
   kGe_{\bc}\cong kNe_{\bc}\otimes_k k^{}_{\alpha}[G/N]
   \]
   where $\alpha\in {\mathrm{H}}^2(G/N,k^{\times})$.  Moreover, the following assertions hold.
   \begin{enumerate}[label={\rm(\alph*)}]
       \item  
       There exists a central extension 
       {$1\longrightarrow Z\longrightarrow H\longrightarrow G/N\longrightarrow 1$}, where $Z$ is the cyclic group of $|G/N|_{p^\prime}$-th roots of unity in $k$, and there exists an idempotent $e\in kZ$ such that $k^{}_{\alpha}[G/N]\cong kHe$.
       \item The block $\bB$ is Morita equivalent to a block $\widehat{\bB}$ of $kH$ having a defect group $\widehat D$ isomorphic to~$D$. After identifying these defect groups, the Morita equivalence is induced by a bimodule with vertex $\Delta D$ and an endo-permutation source.
       \item  If $\alpha$ is trivial, then $H$ can be replaced by $G/N$ itself. 
       If $G/N$ is perfect, then $H$ can be replaced by $[H,H]$ and $Z$ by $Z([H,H])$. 
   \end{enumerate}   
\end{thm}

\noindent Detailed proofs of the main claim, of assertion (a) and of the second claim of assertion  (c)  can be found in \cite[\S3 and \S4]{SalminenSOURCES08}.
    The first claim of assertion (c) follows from \cite[Proposition 1.2.6]{LinckelmannVol1}.
   Assertion~{\rm(b)}, and in particular the more precise information on the vertex and source of the bimodule inducing the Morita equivalence, will be crucial in our subsequent arguments. We record this information in the
following remark. Its proof, together with a proof of Assertion~{\rm(b)}, is given in
Appendix~\ref{app:A}.

\begin{rem}\label{rem:ProofFong2ndRed}%
       Since $\bc$ has defect zero, the block algebra $S:=kNe_{\bc}$ is a matrix algebra. Thus, $S\cong\End_k(V)$ where $V$ is, up to isomorphism, the unique simple $\bc$-module. The module $V$ has a unique extension, up to isomorphism, to a simple $k[ND]$-module, which we also denote by $V$. This extension has vertex $D$, and $\Res_D^{ND}(V)$ is a capped endo-permutation $kD$-module. Moreover, viewed as a $k[\Delta D]$-module via the identification $D\cong\Delta D$, the cap of $\Res_D^{ND}(V)$ is a source of the bimodule inducing the Morita equivalence constructed in the proof of Assertion~{\rm(b)}.
\end{rem}

\vspace{2mm}
\subsection{Splendid stable equivalences of Morita type and trivial source modules}
The following lemma generalizes parts of~\cite[Lemma A.3]{KMN11}. 

\begin{lem}\label{lem:stablesplendid_and_pperm_modules}
     Let $G$ and $H$ be two subgroups of a finite group $L$. Assume  $G$ and $H$ have a common Sylow $p$-subgroup $P$ and the same fusion system on $P$, i.e. $\mathcal{F}_P(G)=\mathcal{F}_P(H)$.
     Let $\bA\in\Bl(kG)$ and let $\bB\in\Bl(kH)$ be blocks with defect group $P$. Assume
     $M$ (resp. $M^{\prime}$) is an indecomposable  finitely generated $(\bA,\bB)$-bimodule (resp. a $(\bB,\bA)$-bimodule),  which,  
     viewed as a $k[G\times H]$-module (resp. as a $k[H\times G]$-module), is  a $p$-permutation module with vertex $\Delta P$, and assume that the pair $(M, M^{\prime})$ induces a stable equivalence of Morita type between $\bA\mathrm{-mod}$ and $\bB\mathrm{-mod}$.
     \par
     Assume, furthermore, that $X$ is a non-projective indecomposable $\bB$-module. Then, the following assertions hold.
     \begin{enumerate}[label={\rm(\alph*)}]
         \item There exists a non-projective indecomposable $\bA$-module $Y$, unique up to isomorphism, such that $M\otimes_{\bB}X = Y\oplus (\mathrm{proj})$ as $\bA$-modules. 
         \item If $X$ is a trivial source module, then so is $Y$. Moreover, $X$ and $Y$ have a vertex in common.
     \end{enumerate}
\end{lem}

\begin{proof}
    Assertion~(a) is well-known: it  follows from \cite[Proposition 4.13.2]{LinckelmannVol1}. 
    We prove (b). As $X$ is non-projective, we may choose a vertex $Q$ of $X$ such that $1\neq Q\leq P$ and we have 
    \[    
    X \big| \Ind_Q^H(k)=kH\otimes_{kQ}k\,.
    \]
    Now, by the assumption, we have $M\mid \Ind_{\Delta P}^{G\times H}(k)$, which is $kG\otimes_{kP}kH$ as  $(kG,kH)$-bimodule.
    Therefore, 
    \[
    Y\,\big| \left(M\otimes_{kH}X\right) \big|\,  kG\otimes_{kP}kH \otimes_{kH}kH\otimes_{kQ}k \,\big|\, \Ind_P^G\left(   
                \Res^H_P(\Ind_{Q}^{H}(k))
            \right)
   \]
   and the Mackey Formula yields 
   \[
    Y \,\big|   \bigoplus_{g\in [Q\backslash G/P]}\Ind_{\lconj{g}{Q}\cap P}^{G}(k)\,.
   \]
   Hence, $Y$ is a trivial source module, and it has a vertex $R$ such that  $R\leq\!\lconj{g}{Q}\cap P$ for some $g\in G$.  It follows that $\lconj{g^{-1}}{R}\leq Q\cap \lconj{g^{-1}}{P}\leq Q\leq P$. Since $\mathcal{F}_P(G)=\mathcal{F}_P(H)$, there exists $h\in H$ such that   $\lconj{g^{-1}}{R}=\lconj{h}{R}$, and so $\lconj{h}{R}$ is also a vertex of $Y$ and $R\leq_H Q$. Since the pair $(M, M^{\prime})$ induces a stable equivalence of Morita type, the same argument with $M^{\prime}$ yields $Q\leq_G R$. It follows that $R=_{H}Q$ and $Q=_GR$. Assertion (b) follows. 
\end{proof}

%% file: section_3_blocks_with_dihedral_defect_grps_prelim_res.tex
\section{Blocks with dihedral defect groups: preliminary results}\label{sec:blocksWithDihedralDefectGroups}

To start with, we state two results, which do not require the defect of the blocks considered to be maximal.
First, we restate \cite[Lemma~4.1]{MacGregor22} in terms of splendid Morita equivalences.

\begin{lem}\label{lem:coveringOddIndex}%
    Let $G$ be a finite group and let $N\unlhd G$ be a normal subgroup such that $|G:N|$ is odd.  If $\bb\in\Bl(kN)$ is a block with dihedral defect groups of order $2^n$ $(n\geq 3)$ and $\bB\in \Bl(kG\mid \bb)$, then the $(\bB,\bb)$-bimodule $e_{\bB}kG e_{\bb}$ and its dual realize a splendid Morita equivalence between $\bB$ and $\bb$. (The blocks $\bb$ and $\bB$ are \emph{isomorphic via restriction} in the terminology of \cite[Theorem~4.1 and Definition~4.3]{HidaKoshitani}.)
\end{lem}

\begin{proof}
    First, observe that \cite[Lemma~4.1]{MacGregor22} proves that restriction from $G$ to $N$ induces a height-preserving bijection between $\Irr(\bB)$ and $\Irr(\bb)$. Thus, by \cite[0.1 Th\'{e}or\`{e}me and 0.2 Th\'{e}or\`{e}me]{BroueIsom90} the $(\bB,\bb)$-bimodule $e_{\bB}kG e_{\bb}$ and its dual realize a Morita equivalence between $\bB$ and $\bb$. This equivalence is a splendid Morita equivalence, since $e_{\bB}kG e_{\bb}$ is a $p$-permutation bimodule. 
\end{proof}


 \begin{lem}\label{lem: nilpotentblock}
 Let $G$ be a finite soluble group. 
 Then, any nilpotent  block $\bB\in\Bl(kG)$ with a defect group $D$ isomorphic to a dihedral $2$-group of order $2^n$ $(n\geq 3)$  is splendidly Morita equivalent to $kD_{2^n}$.
 \end{lem}

\begin{proof}
    Let $\bB\in\Bl(kG)$ be a nilpotent block and let $A$ be a source algebra of $\bB$. By Puig's characterisation of the source algebras of nilpotent blocks \cite{PuigNilpotentBlocks88} (see \cite[Theorem 8.11.5]{LinckelmannVol2}), there exists an endo-permutation $kD$-module $V$ such that 
    $A\cong \End_k(V)\otimes_k kD$ as interior $D$-algebras. Moreover, $D$ is a vertex of $V$ and $V$ is a source of the unique simple $\bB$-module. Thus, $V$ is the trivial $kD$-module by Lemma~\ref{lem:sourcesSimpleDihedral}(b)  and the claim follows. 
 \end{proof}


%% file: section_4_reduction_to_GmodO2.tex
\section{Blocks of maximal defect: reduction to \texorpdfstring{$G/O_{2'}(G)$}{G/O2'(G)}}\label{sec:4}%

\noindent First, we use Fong's reductions to reduce the classification of blocks of maximal defect with dihedral defect groups up to splendid Morita equivalence to cyclic $p'$-central extensions of groups $G$ with dihedral Sylow $2$-subgroups and $O_{2^\prime}(G)=1$. 

\begin{lem}\label{lem:Fong2ndredDihedral}%
     Let $G$ be a finite group with dihedral Sylow $2$-subgroups of order $2^n$, where $n\geq 3$. Then any block $\bB\in\Bl^{\max}(kG)$ is splendidly Morita equivalent to a block $\widehat{\bB}\in\Bl_{\max}(k\widehat{L})$, where $L$ is a finite group with  $O_{2'}(L)=1$ and dihedral Sylow $2$-subgroups of order $2^n$,  and $\widehat{L}$ is a central extension of $L$ by a cyclic group of order $|L|_{2'}$. Moreover, if $H^2(L,k^\times)=1$, then $\widehat{L}$ may be replaced by $L$ itself.
\end{lem}

\begin{proof}
We argue by induction on $|G|$. If $|G|=8$, then $G\cong D_8$ and $\bB=kG$.
Taking $L=\widehat{L}=G$ and $\widehat{\bB}=\bB$ proves
the assertion in this case.
Now assume that $|G|>8$ and that the assertion holds for every group satisfying the hypotheses whose order
is strictly smaller than $|G|$.
Set $N:=O_{2'}(G)$,  let $D\in\Syl_2(G)$, let $\bB\in\Bl^{\max}(kG)$. Let $\bc\in\Bl(kN)$ be a block covered by $\bB$  and let $T$ denote its stabilizer in $G$. 
Suppose first that $T<G$.
By the Fong{--}Reynolds Theorem (Theorem~\ref{thm:FongReynolds}) the block $\bB$
is splendidly Morita equivalent to a block $\widetilde{\bB}\in\Bl(kT)$. Moreover, $\widetilde{\bB}$ has a defect group $D'\leq T$
which is $G$-conjugate to $D$.
Since $\bB$ has maximal defect, $D'$ is a dihedral Sylow $2$-subgroup of $T$
of order $2^n$, and $\widetilde{\bB}$ has maximal defect.
The assertion follows by applying the induction
hypothesis to $T$ and $\widetilde{\bB}$, and using the
transitivity of splendid Morita equivalences.

We may therefore assume that $T=G$, so that $\bc$ is
$G$-stable.
Set $L:=G/N$. Then $O_{2'}(L)=1$. Since $D\cap N=1$, the subgroup $DN/N\in\Syl_2(L)$
is isomorphic to $D$. 
 Then, Theorem~\ref{thm:Fong2ndReduction} (Fong's 2nd reduction) yields 
existence of a Morita equivalence between $\bB$ and
a block $\widehat{\bB}\in\Bl_{\max}(k\widehat{L})$,
where $\widehat{L}$ is a central extension of $L$
by a cyclic group of order $|L|_{2'}$.
Moreover, after identifying the defect groups, this
Morita equivalence is induced by a bimodule with
vertex $\Delta D$ and an endo-permutation source.
This source is constructed as described in Remark~\ref{rem:ProofFong2ndRed}: it can be identified with the cap of the endo-permutation $kD$-module $\Res^{ND}_D(V)$, say $V_0$, and we obtain that $V_0\cong k$ by Lemma~\ref{lem:sourcesSimpleDihedral}(b) as it is a source of a simple module for the $2$-nilpotent (hence soluble) group $ND$ and $D$ is dihedral. This proves that the above Morita equivalence is a splendid Morita equivalence.  
The last claim is immediate from Theorem~\ref{thm:Fong2ndReduction}(c). 
\end{proof}

\noindent Next, as a consequence of \cref{lem:coveringOddIndex}, we can describe all maximal defect $2$-blocks with dihedral defect groups of the finite groups $G$ with $O_{2^{\prime}}(G)=1$ up to splendid Morita equivalence. 

\begin{prop}\label{prop:ClassificBlocksDihedralDefectGroupsMaxDefectO_2prime_is_1}%
Let $G$ be a finite group with dihedral Sylow $2$-subgroups of order $2^n$ $(n\geq 3)$. Write $q = r^f$ for a power of a prime number $r$ with $f$ a positive
integer. If $O_{2'}(G)=1$, then any block $\bB\in\Bl^{\mathrm{max}}(kG)$  is splendidly Morita equivalent to $\bB_0(kG)$. 
Moreover, the following assertions hold.
    \begin{enumerate}[label={\rm(\alph*)}]
        \item If $G=D_{2^n}$ with $n\geq 3$, then $\Bl^{\mathrm{max}}(kG)=\{\bB_0(kG)\}$.
        \item If $n=3$ and $G=\fA_7$, then $\Bl^{\mathrm{max}}(kG)=\{\bB_0(kG)\}$.
        \item If $G=\PSL_2(q)\rtimes C_d$ with $(q-1)_2=2^n$ or $(q+1)_2=2^n$ and $d$ is an odd integer such that $d\,|\,f$, 
        then $|\Bl^{\mathrm{max}}(kG)|=d$ and any  $\bB\in\Bl^{\mathrm{max}}(kG)$ is splendidly Morita equivalent to $\bB_0(k\!\PSL_2(q))$.
        \item If $G=\PGL_2(q)\rtimes C_d$ with $(q-1)_2=2^{n-1}$ or $(q+1)_2=2^{n-1}$ and $d$ is an odd integer such that $d\,|\,f$, 
        then $|\Bl^{\mathrm{max}}(kG)|=d$ and any  $\bB\in\Bl^{\mathrm{max}}(kG)$ is splendidly Morita equivalent to $\bB_0(k\!\PGL_2(q))$.
    \end{enumerate}
 In cases {\rm(c)} and {\rm(d)}, any choice of q satisfying the indicated condition yields a representative of the corresponding splendid Morita equivalence class.      
\end{prop}

\begin{proof}
First, it follows from Subsection~\ref{ssec:dihedralSylow2} that $G$ is one of the groups listed in (a) to (d). 
Let $D\in\Syl_2(G)$. 
The number of blocks of maximal defect is easily calculated as follows. 
First, as $D$ is dihedral of order at least $8$, its automorphism group is a $2$-group, so
\[
N_G(D)= DC_G(D) =D\times O_{2^\prime}(N_G(D))   
\text{ and } 
C_G(D)=Z(D)\times O_{2^\prime}({N_G(D)})\,.
\]
Then  \cite[Lemma 2.1]{NavarroTiep12} shows (as a corollary of Brauer's 1st Main Theorem)  that  {$|\Bl^{\mathrm{max}}(kG)|$} is the number of $N_G(D)$-orbits on the conjugacy classes of elements in $O_{2^\prime}({N_G(D)})$.
\par
Now, Case (a) is clear. Indeed, as $G=D$ is a $2$-group, we have $kG=\bB_0(kG)$. Case~(b) is also clear as $D$ is self-normalising in $G=\fA_7$, and so the principal block is the unique block of maximal defect by the above.\par
Thus, we can assume that $G$ is one of the groups in (c) or (d).  Set $G_0:=\PSL_2(q)$ in Case~(c) and $G_0:=\PGL_2(q)$ in Case (d). Then $D\in\Syl_2(G_0)$ and  we have $N_{G_0}(D)=D$. See~\cite[Theorem~1.4.3(d)]{BonnafeBookSL2q} for $\PSL_2(q)$. 
For $G_0=\PGL_2(q)$ this easily follows from the fact that we may choose $D$ to be the Sylow $2$-subgroup of the normalizer of the split (resp. non-split) torus $T$ when $q\equiv 1\pmod{4}$ (resp. $q\equiv 3\pmod{4}$) and $N_{G_0}(T)\cong D_{2(q-1)}$ (resp. $N_{G_0}(T)\cong D_{2(q+1)}$). Any element of $N_{G_0}(D)$ must then normalize the unique cyclic subgroup of $D$ of index $2$ and it follows that $N_{G_0}(D)\leq N_{G_0}(T)$ and $N_{G_0}(D)=D$.
It follows that the principal block is the unique block of maximal defect of $kG_0$. If $d>1$, then $N_G(D)=D\times C_d$ and we obtain that $|\Bl^{\mathrm{max}}(kG)|=d$ by the above counting argument. Moreover, any block $\bB\in\Bl^{\mathrm{max}}(kG)$ 
covers a block $\bb\in\Bl^{\mathrm{max}}(kG_0)=\{\bB_0(kG_0)\}$ (as a defect group of $\bb$ is $G$-conjugate to $D\cap G_0=D$). By Lemma~\ref{lem:coveringOddIndex}
any such block is splendidly Morita equivalent to $\bB_0(kG_0)$. In particular, all the elements of $\Bl^{\mathrm{max}}(kG)$ are splendidly Morita equivalent to $\bB_0(kG)$.
Finally,  in cases {\rm(c)} and {\rm(d)}, the fact that the corresponding splendid Morita equivalence class does not depend on $q$ is given by \cite[Theorem~1.1]{KoLa19}.
The claim follows.
\end{proof}

%% file: sections_5_6_and_7.tex
\section{The faithful blocks of maximal defect of \texorpdfstring{$3.\fA_7$}{3.A7}}\label{sec:3A7}%

In this section and in the next section, we prove that any block of maximal defect arising from the non-split central extensions originating from the non-trivial Schur multipliers in Lemma~\ref{lem:SchurMult_2prime} is again splendidly Morita equivalent to a principal block. Our arguments require the known character and block theory of these groups, for which we use the notation of \cite{ATLAS} and \cite{ModularAtlas}. 
For the purpose of our calculations in this section, we assume that $(K,\mathcal{O},k)$ is a $2$-modular system which is large enough for $3.\fA_7$ and all of its subgroups (this includes the groups $3.\fA_6$ and $\PSL_2(7)\cong \GL_3(2)$, which appear in Section~\ref{sec:3A6})   
and where $k$ is our previously fixed algebraically closed field of characteristic~$2$.\\

\noindent The triple cover $3.\fA_7$ of $\fA_7$ possesses two faithful $2$-blocks of maximal defect, denoted Block 3 and Block 4 in \cite[A7mod2]{ModularAtlas}, which are contragredient to each other. In the following, we prove that these blocks are splendidly Morita equivalent to $k\fS_4$.
Recall here that $\fS_4\cong\PGL_2(3)$. 
It is clear from Erdmann's work \cite{ErdmannTameBook} on tame blocks together with \cite[Theorem~2.1]{MacGregor22} that these blocks are Morita equivalent, because they have the same decomposition matrix.  Erdmann's approach, classifying basic algebras, however, does not provide us with a bimodule realising the equivalence. We construct such a bimodule using the Brauer correspondence.

\vspace{2mm}
\subsection{Context and notation}
We first introduce the notation necessary to solve the problem. 

\begin{nota}\label{nota:Section_3A7}%
Throughout this subsection, we use the following notation. First, 
\begin{enumerate}[label=(\arabic*)]
    \item $G:=3.\fA_7$ denotes the triple cover of $\fA_7$;
    \item $Z:=Z(G)=:\langle z\rangle\cong C_3$\,;
    \item  $P\in\Syl_2(G)$ is a fixed Sylow $2$-subgroup,  
      so $P\cong D_8$ and $N_G(P)=Z\times P\cong C_3\times D_8$\,;
    \item  $Q\lneq P$ denotes the unique subgroup such that $Q\cong C_2\times C_2$ and $L:=N_G(Q)=Z\times H$ with $H\cong \fS_4$, so $L\cong C_3\times \fS_4$;  
    \item $Q_2\leq P$ denotes the second subgroup of $P$ such that $Q_2\cong C_2\times C_2$, and its normaliser is $N_G(Q_2)\cong ((C_2\times C_2):( (C_3\times C_3):C_3 )):C_2$.
\end{enumerate}
Thus, we have the following chain of subgroups:
 $$Q\leq P\leq N_G(P)=Z\times P\leq L=Z\times H\leq G\,.$$
Next, we consider the following blocks: 
\begin{enumerate}[label=(\arabic*), start=6]
    \item $\bA:=$ Block 3 of $kG$ according to the notation of \cite[A7mod2]{ModularAtlas},\\
    $\bA^{\ast}:=$ Block 4 of $kG$ according to the notation of \cite[A7mod2]{ModularAtlas};   
    \item\label{nota:3A7_B0(kH)=kH}  $\bB_0(kH)=kH$, the principal block of $H\cong \fS_4$ (all partitions of $4$ have the same $2$-core, namely the empty set, so it is clear that $kH$ has just one block);
    \item\label{nota:3A7_def_c_c3} $\bc$ and $\bc^\ast$,  the Brauer correspondents in $L=N_G(Q)$ of $\bA$ and $\bA^\ast$ respectively.
\end{enumerate}
Moreover,  write:
\begin{enumerate}[label=(\arabic*), start=9]
    \item\label{nota:blocksOfkZ}  $\Irr(Z):=\{\zeta_0:=1_{C_3},\zeta_1,\zeta_2\}$ and for each $0\leq i\leq 2$ let $e_i:=e_{\zeta_i}\in kZ$ be the block idempotent associated to $\zeta_i$. (See \cite[Chapter~3, Theorem~6.22]{NagaoTsushimaBook}.) Therefore, $kZ=kZe_0\oplus kZe_1\oplus kZe_2$ is the block decomposition of $kZ$. 
 \end{enumerate} 
Clearly, the Sylow $2$-subgroup $P$ is a defect group of the blocks $\bA$, $\bA^\ast$, $\bB_0(kH)$, $\bc$ and $\bc^\ast$, whereas the blocks of $kZ$ are of  defect zero.\\ 
The simple $kH$-modules and the simple $\bA$-modules are as follows.
\begin{enumerate}[label=(\arabic*),start=10]
    \item\label{nota:simplesInkH} We can read from~\cite[Appendix I.E]{JamesKerberBook} and/or  the Character Table Library package \verb|CTblLib| in GAP4~\cite{GAP} that the block $\bB_0(kH)=kH$ has, up to isomorphism, two simple modules; the trivial module and a simple module of dimension $2$, which are both liftable. Hence, we set $\Irr_k(\bB_0(kH))=:\{k_H,2_H\}$.
    \item\label{nota:3A72-decompMatrixA} By~\cite[A7mod2, Block~3]{ModularAtlas} the blocks $\bA$ and $\bA^\ast$ both have two simple modules: one of dimension 15 and one of dimension 6. Again, identifying these with their dimensions, we set $\Irr_k(\bA):=\{15, 6\}$  and $\Irr_k(\bA^\ast):=\{15^\ast, 6^\ast\}$. The $2$-decomposition matrix of $\bA$ and the values at the conjugacy classes of $2$-elements are: 
        \[
        \begin{tabular}{c:c|cc}{\footnotesize{\cite[p.10]{ATLAS}}} &  \cite{ModularAtlas} & $15$ & $6$\\
            \hline  
            $\chi^{}_{17}$ & $\chi^{}_{6_2}$  & $0$ & $1$\\
            $\chi^{}_{18}$ & $\chi^{}_{15_2}$ & $1$ & $0$\\
            $\chi^{}_{19}$ & $\chi^{}_{15_4}$ & $1$ & $0$\\
            $\chi^{}_{20}$ & $\chi^{}_{21_2}$ & $1$ & $1$\\
            $\chi^{}_{21}$ & $\chi^{}_{21_4}$ & $1$ & $1$\\
        \end{tabular}
        \qquad \qquad
        \begin{tabular}{c|ccc} &  $1A$ & $2A$ & $4A$\\
            \hline  
            $\chi^{}_{17}$ & $6$  & $2$ & $0$\\
            $\chi^{}_{18}$ & $15$ & $-1$ & $-1$\\
            $\chi^{}_{19}$ & $15$ & $3$ & $1$\\
            $\chi^{}_{20}$ & $21$ & $1$ & $-1$\\
            $\chi^{}_{21}$ & $21$ & $-3$ & $1$\\
        \end{tabular}
        \]
\end{enumerate} 
 The block $\bA$, seen as a left $k[G\times G]$-module, is indecomposable with vertex $\Delta P$ and trivial source. (See e.g.~\cite[Chap. 5, Theorem~10.8]{NagaoTsushimaBook}.) Moreover, by the definitions we have $N_{G\times G}(\Delta P)\leq N_G(P)\times N_G(P)$, so 
    \(
        N_{G\times G}(\Delta P)\leq G\times L
    \)    
        and
    \(  
    N_{G\times G}(\Delta P)\leq L\times G
    \)\,.
    Therefore, we may consider the Green correspondence $\mathfrak{f}:=\mathfrak{f}(G\times G,\Delta P,L\times G)$ (following the notation of \cite[Chap. 4, \S4]{NagaoTsushimaBook}). Thus, we may set 
\begin{enumerate}[label=(\arabic*),start=12]
    \item\label{nota:3A7_candidate_bimodule} $M:=\mathfrak{f}(\bA)$, i.e. the unique indecomposable direct summand of $\Res^{G\times G}_{L\times G}(\bA)$ with vertex~$\Delta P$, which is a $(\bc,\bA)$-bimodule.
\end{enumerate}   
\end{nota}

\vspace{2mm}
\subsection{A stable
equivalence of Morita type}\label{ssec:stEqMorTp}
We first prove that the $(\bc,\bA)$-bimodule $M$ (and its dual) realise a stable equivalence of Morita type between $\bA$ and $\bc$ using Linckelmann's glueing method \cite[Theorem~3.1]{Lin01}. To this end, we collect information about the structure of the centralizers of the $2$-subgroups below and in the next remark, which can be found at \cite[Groups/Abstract/7560.g]{lmfdb} or using GAP4~\cite{GAP}.
\begin{nota}\label{nota:Section_3A7_part_II}\leavevmode
\begin{enumerate}[label=(\arabic*),start=13]
    \item\label{nota:local-3A7} 
    We have 
    \(
C_L(P)=Z\times Z(P)=C_G(P)
\)
and, under the identification
\[
k[Z\times Z(P)]\cong kZ\otimes_k kZ(P)\,,
\]
we have
\[
kC_L(P)\Br_P(e_\bc)=k[Z\times Z(P)](e_1\otimes 1_{kZ(P)})=kC_G(P)\Br_P(e_\bA)\,.
\]
We identify $e_1$ with $e_1\otimes 1_{kZ(P)}$. Hence,
\[
f_P:=\Br_P(e_\bc)=e_1=\Br_P(e_\bA)=: e_P\,.
\]
In particular, $(P,f_P)$ and $(P,e_P)$ are maximal Brauer pairs of $\bc$ and $\bA$, respectively.
\item For any $R\le P$, let $f_R$ (resp. $e_R$) denote the unique block idempotent of $kC_L(R)$ (resp. $kC_G(R)$) such that 
\[(R,f_R)\le_L(P,f_P)\qquad\text{and}\qquad (R,e_R)\le_G (P,e_P)\,.
\] 
The fusion systems of $\bc$ and $\bA$ with respect to the maximal Brauer pairs $(P,f_P)$ and $(P,e_P)$ are isomorphic, as both blocks have the same number of isomorphism classes of simple modules. 
(This follows from~\cite[Part IV, Section 5.7]{AschKessOli11} together
with~\cite{BrauerDihedral}; see also \cite[Theorem 8.1(2)]{Sambale14}.)
We identify them and denote the resulting fusion system on $P$ by~$\calF$. The group $P$ has two $\calF$-conjugacy classes of involutions and
\[
\Aut_\calF(Q)\cong \fS_3\quad\text{and}\quad \Aut_\calF(Q_2)\cong C_2\,.
\]
\end{enumerate}
\end{nota}

\begin{rem}\label{rem:local-3A7}
    We record some local information that will be used repeatedly in the proof of Proposition~\ref{prop:gluing-3A7} and  which can be read off from \cite[Groups/Abstract/7560.g]{lmfdb}.  
    \begin{enumerate}[label=(\alph*)]
    \item
    Let $R\le Q_2$ be a subgroup of order $2$ such that $R\neq Z(P)$. Then, $C_G(R)=N\rtimes \hat{P}$
    for some $\hat{P}\in \Syl_2(G)\setminus \{P\}$ containing $Q_2$, where 
    \(
    N:=O_3(C_G(R))\cong C_3\times C_3\,.
    \)
    As $C_G(R)$ contains $C_G(Q_2)$, we obtain that
    \[
    C_G(Q_2)=C_{C_G(R)}(Q_2)=C_{N\rtimes\hat{P}}(Q_2)=C_N(Q_2)\rtimes C_{\hat{P}}(Q_2)= C_N(Q_2)\rtimes Q_2\,.
    \]
    Since $|C_G(Q_2)|=36$, we get that $C_N(Q_2)=N$, and hence $C_G(Q_2)=N\times Q_2$.
    
    Next, as $Z(P)\leq Q_2$, we have $N\leq C_G(Q_2)\leq C_G(Z(P))$
    and clearly $P\leq C_G(Z(P))$. Moreover, $Q_2\unlhd P$, so $P$ normalizes $C_G(Q_2)$ and therefore also its characteristic subgroup $N=O_3(C_G(Q_2))$. Thus, $N\rtimes P\leq C_G(Z(P))$. Since $|C_G(Z(P))|=72=|N\rtimes P|$, we conclude that $C_G(Z(P))=N\rtimes P$. 

    \item
    Next, we record the corresponding local block idempotents. For every non-trivial subgroup $S\le P$, one has $C_L(S)=Z\times C_P(S)$. Since $C_P(S)$ is a $2$-group, under the natural identification $kC_L(S)\cong kZ\otimes_k kC_P(S)$ we have $f_S=e_1\otimes 1_{kC_P(S)}$. As before, we identify $e_1$ with
    $e_1\otimes 1_{kC_P(S)}$ and simply write $f_S=e_1$.
    
    Since $C_G(Q_2)=N\times Q_2$ and $Q_2$ is a $2$-group, there is a unique $\varepsilon\in\Bl(kN)$ such that $e_{Q_2}=\varepsilon\otimes 1_{kQ_2}$. Identifying $\varepsilon$ with $\varepsilon\otimes 1_{kQ_2}$ as above, we simply write $e_{Q_2}=\varepsilon$. Since $e_{Q_2}$ lies over the block $e_1$ of $kZ$, we have $\varepsilon e_1=\varepsilon$. 
    
    Now, the block $e_{Z(P)}$ of $kC_G(Z(P))$ covers the same block $\varepsilon$ of $kN$. We have $C_G(Z(P))=N\rtimes P$, and $e_{Z(P)}$ has defect group $C_P(Z(P))=P$. It follows that $\varepsilon$ is $P$-stable, and hence $e_{Z(P)}=\varepsilon$.
    
    Finally, let $R\leq Q_2$ be a subgroup of order $2$ such that $R\neq Z(P)$. The block $e_R$ of $kC_G(R)$ also covers $\varepsilon$. Recall from part~{\rm (a)} that $C_G(R)=N\rtimes\hat P$
    for some $\hat P\in\Syl_2(G)$ containing $Q_2$. Since $e_R$ has defect group $C_P(R)=Q_2$, the inertia group of $\varepsilon$ in $\hat P$ is $Q_2$. Thus $\varepsilon$ is not $\hat P$-stable and, for any $x\in\hat P\setminus Q_2$, we have $e_R = \operatorname{Tr}_{Q_2}^{\hat P}\varepsilon=\varepsilon+{}^x\varepsilon$. 
\end{enumerate}    
\end{rem}

\begin{rem}
Notice that for every $R\leq P$, we have
\[
C_L(R)\times C_G(R)=C_{L\times G}(\Delta R)
\leq 
N_{L\times G}(\Delta R)
\leq N_L(R)\times N_G(R)\,.
\]
If $R\cong C_2$, then any element normalising $R$ also centralises $R$, thus in this case 
\[
C_{L\times G}(\Delta R)
=C_L(R)\times C_G(R) 
= N_L(R)\times N_G(R) 
=N_{L\times G}(\Delta R)\,.
\]
It follows that if $Y$ is a $k[L\times G]$-module, its Brauer construction $Y[\Delta R]$ can be considered as a $k[C_L(R)\times C_G(R)]$-module via restriction, hence as a $(kC_L(R),kC_G(R))$-bimodule.
We use this remark repeatedly in the proof of Proposition~\ref{prop:gluing-3A7} below without further reference.
\end{rem}

\begin{prop}\label{prop:gluing-3A7}
The $(\bc,\bA)$-bimodule $M$  defined in Notation~\ref{nota:Section_3A7}{\rm\ref{nota:3A7_candidate_bimodule}} induces a stable equivalence of Morita type between the blocks $\bA$ and $\bc$. 
\end{prop}
\begin{proof}
We prove the proposition using the glueing method given by~\cite[Theorem~3.1]{Lin01}. Since the blocks $\bc$ and $\bA$ have isomorphic fusion systems, it suffices to show that for any nontrivial subgroup $R\le P$, the bimodule $f_RM[\Delta R]e_R$ induces a splendid Morita equivalence between $kC_L(R)f_R$ and $kC_G(R)e_R$.
\par
To this end, consider the $(\bc,\bA)$-bimodule $\widetilde{M}:=e_\bc kGe_\bA$. We have $\widetilde{M} \mid \Res^{G\times G}_{L\times G}(kG e_\bA)$. Moreover, the Brauer construction
\[
\widetilde{M}[\Delta P]\cong \Br_P(e_\bc)kC_G(P)\Br_P(e_\bA)=k[Z\times Z(P)]e_1
\]
is non-zero. It follows that $\widetilde{M}$ has an indecomposable direct summand with vertex $\Delta P$. The uniqueness of the Green correspondent implies that $M \mid \widetilde{M}$. Moreover, since $M$ is the unique indecomposable summand of $\widetilde{M}$ with a vertex $\Delta P$, the bimodule 
\[
f_PM[\Delta P]e_P=\Br_P(e_\bc)M[\Delta P]\Br_P(e_\bA)=\Br_P(e_\bc)\widetilde{M}[\Delta P]\Br_P(e_\bA)
\]
is non-zero. It follows that for any $R\le P$, the bimodule $f_RM[\Delta R]e_R$ is a direct summand of $f_R\widetilde{M}[\Delta R]e_R$. Moreover, this summand is non-zero. Indeed, by \cite[Proposition~5.3(c)]{BoltjePerepelitsky}, the pair $(\Delta P, f_P\otimes e_P^*)$ is a maximal $M$-Brauer pair. By \cite[Remark~10.2(c)]{BoltjePerepelitsky}, we have 
\[
(\Delta R, f_R\otimes e_R^*)\le (\Delta P,f_P\otimes e_P^*)\,.
\]
Hence by \cite[Proposition~5.3(b)]{BoltjePerepelitsky}, the pair $(\Delta R, f_R\otimes e_R^*)$ is an $M$-Brauer pair as well. Thus, $f_RM[\Delta R]e_R$ is non-zero. We now verify the required local splendid Morita equivalences case by case.

\medskip
\noindent{{\bf Case \rm (i).}} Suppose that $R$ is one of $P$, $Q$ or the cyclic subgroup of order $4$. Then in each of these cases, we have 
\[
C_L(R)=Z\times Z(R)=C_G(R)\,.
\]
Indeed, $L=N_G(Q)$ implies that $C_G(Q)=C_L(Q)=Z\times Q$, and $N_G(P)=Z\times P$ implies that $C_G(P)=Z\times Z(P)=C_L(P)$. Moreover, if $R\cong C_4$, then $C_{\fA_7}(\overline{R})=\overline{R}$ implies that $C_G(R)=Z\times R=C_L(R)$. Now, we have
\[
kC_L(R)f_R=k[Z\times Z(R)]e_1=kC_G(R)e_R
\]
and hence the bimodule
\[
f_R\widetilde{M}[\Delta R]e_R=e_1k[Z\times Z(R)]e_1=k[Z\times Z(R)]e_1
\]
induces a splendid Morita equivalence between $kC_L(R)f_R$ and $kC_G(R)e_R$. In particular, the bimodule $f_R\widetilde{M}[\Delta R]e_R$ is indecomposable and thus $f_RM[\Delta R]e_R=f_R\widetilde{M}[\Delta R]e_R$.

\medskip
\noindent{{\bf Case \rm (ii).}} Suppose that $R=Q_2$. By Remark~\ref{rem:local-3A7},
\[
C_L(R)=Z\times R,\qquad C_G(R)=N\times R,
\qquad f_R=e_1,\qquad e_R=\varepsilon.
\]
Therefore, the bimodule
\[
f_R\widetilde{M}[\Delta R]e_R=e_1k[NR]\varepsilon
\]
induces a splendid Morita equivalence between $kC_L(R)f_R=k[ZR]e_1$ and $kC_G(R)e_R=k[NR]\varepsilon$. As in Case~{\rm (i)}, it follows that the bimodule $f_R\widetilde{M}[\Delta R]e_R$ is indecomposable. Thus we have $f_RM[\Delta R]e_R=f_R\widetilde{M}[\Delta R]e_R$.

\medskip
\noindent{{\bf Case \rm (iii).}} Suppose that $R=Z(P)$. By Remark~\ref{rem:local-3A7},
\[
C_L(R)=Z\times P,\qquad C_G(R)=N\rtimes P,
\qquad f_R=e_1,\qquad e_R=\varepsilon.
\]
Therefore,
\[
f_R\widetilde{M}[\Delta R]e_R=e_1k[N\rtimes P]\varepsilon
\]
induces a splendid Morita equivalence between
$k[Z\times P]e_1$ and $k[N\rtimes P]\varepsilon$.
Hence it is indecomposable, and
\[
f_RM[\Delta R]e_R
   =f_R\widetilde{M}[\Delta R]e_R.
\]

\medskip
\noindent{{\bf Case \rm (iv).}} Finally, suppose that $R\cong C_2$ and that $R$ is  not $\calF$-conjugate to $Z(P)$. By Remark~\ref{rem:local-3A7},
\[
C_L(R)=Z\times Q_2\,,\quad C_G(R)=N\rtimes \hat{P},\qquad f_R=e_1\,,
\]
where $\hat{P}\in \Syl_2(G)$ contains $Q_2$, and $e_R=\varepsilon+{}^{x}\varepsilon$ for any $x\in \hat{P}\setminus Q_2$. Recall also that $e_{Q_2}=\varepsilon$ and $e_1\varepsilon=\varepsilon$. 

By the Fong-Reynolds Theorem, the bimodule $\varepsilon k[C_G(R)]e_R$ induces a splendid Morita equivalence between $k[NQ_2]\varepsilon$ and $kC_G(R)e_R$. Moreover, the bimodule $e_1k[NQ_2]\varepsilon$ induces a splendid Morita equivalence between $k[ZQ_2]e_1$ and $k[NQ_2]\varepsilon$. Hence their composition
\[
e_1k[NQ_2]\varepsilon\otimes_{k[NQ_2]\varepsilon}\varepsilon kC_G(R)e_R\cong \varepsilon kC_G(R)e_R
\]
induces a splendid Morita equivalence between the blocks $k[ZQ_2]e_1=kC_L(R)f_R$ and $k[N\hat{P}]e_R=kC_G(R)e_R$. In particular, the $(kC_L(R)f_R, kC_G(R)e_R)$-bimodule $\varepsilon kC_G(R)e_R$ is indecomposable. Note that
\[
e_1e_R=e_1(\varepsilon+{}^{x}\varepsilon)=e_1\varepsilon+e_1{}^{x}\varepsilon=\varepsilon+{}^{x}\varepsilon=e_R
\]
and hence
\[
f_R\widetilde{M}[\Delta R]e_R=e_1kC_G(R)e_R=e_1e_RkC_G(R)e_R=\varepsilon kC_G(R)e_R\oplus {}^{x}\varepsilon kC_G(R)e_R
\]
is a sum of two indecomposable modules. We need to show that $f_RM[\Delta R]e_R\cong \varepsilon kC_G(R)e_R$.

Set $Y:={}^{x}\varepsilon kC_G(R)e_R$ viewed as a $(kC_L(R),kC_G(R))$-bimodule. Since $f_RM[\Delta R]e_R$ is non-zero, it is enough to prove that $Y$ is not a direct summand of $f_RM[\Delta R]e_R$. Suppose, for a contradiction, that it is. Note that $\varepsilon kC_G(R)e_R$ has vertex $\Delta Q_2$ and ${}^{x}\varepsilon kC_G(R)e_R$ has vertex $\Delta_\varphi(Q_2)=\Delta(Q_2, \varphi, Q_2)$ where $\varphi=i_x\mid_{Q_2}:Q_2\longrightarrow Q_2$. 

Since $R\cong C_2$, we have $N_{L\times G}(\Delta R)=C_L(R)\times C_G(R)$, and by \cite[Proposition~5.8.7]{LinckelmannVol1},
\[
M[\Delta R]\mid
\Res_{C_L(R)\times C_G(R)}^{L\times G}(M).
\]
Hence
\[
Y\mid
\Res_{C_L(R)\times C_G(R)}^{L\times G}(M).
\]
The Brauer construction
$Y[\Delta_\varphi(Q_2)]$ is naturally a module for $N_{C_L(R)\times C_G(R)}(\Delta_\varphi(Q_2))$, via inflation. Restricting to the centralizer, and using
\[
C_{C_L(R)\times C_G(R)}(\Delta_\varphi(Q_2))=C_L(Q_2)\times C_G(Q_2),
\]
we see that $f_{Q_2}\otimes e_{Q_2}^*$ acts naturally on $Y[\Delta_\varphi(Q_2)]$. We claim that $(f_{Q_2}\otimes e_{Q_2}^*)Y[\Delta_\varphi(Q_2)]\neq 0$. In fact, since ${}^{x}\varepsilon\otimes e_R^*$ is fixed by $C_L(R)\times C_G(R)$, by \cite[Lemma~3.7]{BoltjePerepelitsky}, we have an isomorphism
\[
Y[\Delta_\varphi(Q_2)]\cong \Br_{\Delta_\varphi(Q_2)}({}^{x}\varepsilon\otimes e_R^*) kC_G(R)[\Delta_\varphi(Q_2)]
\]
of $kN_{C_L(R)\times C_G(R)}(\Delta_\varphi(Q_2))$-modules. Note that
\[
N_{C_L(R)\times C_G(R)}(\Delta_\varphi(Q_2))=C_L(R)\times C_G(Q_2).
\]
Now both ${}^{x}\varepsilon$ and $e_R$ are in $kN$, and $Q_2$ centralizes $N$. Hence $\Br_{\Delta_\varphi(Q_2)}
({}^{x}\varepsilon\otimes e_R^*)={}^{x}\varepsilon\otimes e_R^*$. Moreover, the $\Delta_\varphi(Q_2)$-fixed points of $C_G(R)$ are precisely $xC_G(Q_2)$. Indeed, for $g\in C_G(R)$ and $u\in Q_2$, one has $(\varphi(u),u)g=g$ if and only if $xux^{-1}gu^{-1}=g$. This holds for all $u\in Q_2$ if and only if $x^{-1}g\in C_G(Q_2)$. Thus, $kC_G(R)[\Delta_\varphi(Q_2)] \cong k[xC_G(Q_2)]$. Therefore,
\begin{align*}
	Y[\Delta_\varphi(Q_2)]&\cong {}^{x}\varepsilon k[xC_G(Q_2)]e_R=x\varepsilon kC_G(Q_2)e_R=x\varepsilon kC_G(Q_2)
\end{align*}
where the last equality follows from
$\varepsilon e_R=\varepsilon$. Finally, since $f_{Q_2}=e_1$, $e_{Q_2}=\varepsilon$, and
$e_1\varepsilon=\varepsilon$, we obtain
\[
(f_{Q_2}\otimes e_{Q_2}^*)
Y[\Delta_\varphi(Q_2)]\cong x\varepsilon kC_G(Q_2)\neq0.
\]
Since $Y$ is a direct summand of
$\Res_{C_L(R)\times C_G(R)}^{L\times G}(M)$, it follows that
\[
(f_{Q_2}\otimes e_{Q_2}^*) M[\Delta_\varphi(Q_2)]\neq0.
\]
Moreover, by \cite[Remark~10.2(a)]{BoltjePerepelitsky}, $\bigl(\Delta_\varphi(Q_2),f_{Q_2}\otimes e_{Q_2}^*\bigr)$ is a Brauer pair of $L\times G$, and so by definition, it is an $M$-Brauer pair. By \cite[Proposition~5.3(c)]{BoltjePerepelitsky}, this $M$-Brauer pair is contained in a maximal $M$-Brauer pair. Since $(\Delta P,f_P\otimes e_P^*)$
is a maximal $M$-Brauer pair and all maximal $M$-Brauer pairs are $(L\times G)$-conjugate, there exist $l\in L$ and $g\in G$ such that
\[
\bigl(\Delta_\varphi(Q_2),f_{Q_2}\otimes e_{Q_2}^*\bigr)\leq{}^{(l,g)}\bigl(\Delta P,f_P\otimes e_P^*\bigr).
\]
It follows that there exists a subgroup $U\leq P$ such that
\[
\Delta_\varphi(Q_2)={}^{(l,g)}\Delta U\quad \text{and} \quad f_{Q_2}\otimes e_{Q_2}^*={}^{(l,g)}(f_U\otimes e_U^*)={}^{l}f_U\otimes {}^{g}e_U^*\,.
\]
This means that
	\[
	{}^lU=Q_2={}^gU,\qquad f_{Q_2}={}^{l}f_U,\qquad e_{Q_2}={}^{g}e_U\qquad  \text{and}\qquad \varphi=i_l\circ i_g^{-1}:Q_2\longrightarrow Q_2\,.
	\]
Therefore, the maps $i_l|_U: U\longrightarrow Q_2$ and $i_g|_U: U\longrightarrow Q_2$ are isomorphisms in the fusion systems of $\bc$ and $\bA$, respectively. Under our identification of these two fusion
systems, both are isomorphisms in $\calF$. Since $Q$ and $Q_2$ are not $\calF$-conjugate, it follows that $U=Q_2$. Hence both $i_l|_{Q_2}$ and $i_g|_{Q_2}$ are in $\Aut_\calF(Q_2)$,
and therefore $\varphi\in\Aut_\calF(Q_2)$. But $\Aut_\calF(Q_2)=\Aut_P(Q_2)$.

Now $\varphi$ is non-trivial since $\varphi=i_x|_{Q_2}$ for $x\in \hat{P}\setminus Q_2$. On the other hand, $\hat P\leq C_G(R)$, and hence $x$ centralizes $R$. Thus $\varphi=i_x|_{Q_2}$ fixes $R$. This gives a contradiction. Indeed, the unique non-trivial element of $\Aut_P(Q_2)$ fixes $Z(P)$ and interchanges the other two subgroups of order $2$ in $Q_2$. Since $R\neq Z(P)$, no non-trivial element of $\Aut_P(Q_2)=\Aut_\calF(Q_2)$ can fix $R$.

Therefore, $Y$ is not a direct summand of
$f_RM[\Delta R]e_R$ and we conclude that
\[
f_RM[\Delta R]e_R\cong \varepsilon kC_G(R)e_R.
\]
Hence $f_RM[\Delta R]e_R$ induces the required splendid Morita equivalence.
\end{proof}

\vspace{2mm}
\subsection{The splendid Morita equivalences}
It now remains to prove that the stable equivalence of Proposition~\ref{prop:gluing-3A7} is a splendid Morita equivalence. We achieve this using a method due to Linckelmann (see \cite[Theorem 2.1(ii) and~(iii)]{Lin96StableEqs}), which consists in proving that the functor $M\otimes_{\bA}-$ maps the simple $\bA$-modules to the simple $\bc$-modules. Thus, we first need to collect information on the simple modules of the blocks involved.

\begin{lem}\label{lem:SplendidlyEquivBrauerCorrespondents}%
The following assertions hold.
\begin{enumerate}[label={\rm(\alph*)}]
    \item
We have  $|\Bl(kL)|=3$. Furthermore, 
    \[    \bB_0(kL)=kZe_0\otimes_k\bB_0(kH),\quad \bc=kZe_1\otimes_k\bB_0(kH),\text{ and } \bc^\ast=kZe_2\otimes_k\bB_0(kH)
    \]
    (permuting the labelling of $e_1$ and $e_2$ if necessary)  
    and all three blocks are splendidly Morita equivalent to $\bB_0(kH)$.   
    \item The block $\bc$ contains precisely two pairwise non-isomorphic simple modules, one of dimension $1$ and one of dimension $2$. We identify these modules with their dimensions and write $\Irr_k(\bc)=:\{1,2\}$.
    \item 
    The simple $\bc$-module $1$ is the unique trivial source $\bc$-module with vertex $P$.
    The simple $\bc$-module~$2$ is a trivial source module with vertex $Q$. 
\end{enumerate}    
\end{lem}

\begin{proof}
\begin{enumerate}[label={\rm(\alph*)}]
    \item Since $Z$ is a $2'$-group and  $|\Irr(Z)|=3$,  we have  $|\Bl(kZ)|=3$ and each block of $kZ$ is of the form $kZe_i$ for some $i\in\{0,1,2\}$ (as described in Notation~\ref{nota:Section_3A7}\ref{nota:blocksOfkZ}). By~Notation~\ref{nota:Section_3A7}\ref{nota:3A7_B0(kH)=kH}, $\Bl(kH)=\{\bB_0(kH)\}$. Thus, it is clear from Subsection~\ref{ssec:blocks_of_central_producrs} that $|\Bl(kL)|=3$ and any block of $kL$ has the given form. The last claim follows from Lemma~\ref{lem:SplendidMoritaEquivalence_direct_and_central_product}.
    \item By (a), any simple $\bc$-module is of the form $S_1\otimes_kT$, where $S_1$ is the unique simple $kZe_1$-module and $T$ is a simple $\bB_0(kH)$-module. Now, $\dim_k(S_1)=1$ and by Notation~\ref{nota:Section_3A7}\ref{nota:simplesInkH} we have $\Irr_k(\bB_0(kH))=\{k_H,2_H\}$. Hence $1:=S_1\otimes_k k_H$ is a simple module of dimension~$1$ and $2:=S_1\otimes_k 2_H$ is a simple module of dimension $2$. The claim follows.
    \item It is clear that $k_H$ is a trivial source $kH$-module with vertex~$P$. 
    Next, as $H\cong \fS_4$, $Q\unlhd H$ and   $\fS_4/Q\cong \fS_3$, any trivial source $kH$-module with vertex $Q$ is the inflation from $H/Q$ to $H$ of a PIM (see e.g. \cite[Theorem~4.6(c)]{LasSurvey}). Moreover, we can read from \cite[Appendix I.E]{JamesKerberBook} and/or the Character Table Library package \verb|CTblLib| in GAP4~\cite{GAP} that  the $2$-decomposition matrix of $k\fS_3$ is 
    \[
        \begin{tabular}{c:c|cc}{\footnotesize{partition notation}} &  {\footnotesize{degree}} & $k_{\fS_3}$ & $2_{\fS_3}$\\
            \hline  
            $\chi^{}_{(3)}=1_{\fS_3}$  & $1$   & $1$ & $0$\\
            $\chi^{}_{(1^3)}=: \mathrm{sgn}_{\fS_3}\text{\;(sign character)}$& $1$& $1$ & $0$\\
            $\chi^{}_{(2,1)}$& $2$ & $0$ & $1$\\
        \end{tabular}
    \]
    where $2_{\fS_3}$ denotes a simple $k\fS_3$-module of dimension $2$. 
    It follows that $2_H=\Inf_{H/Q}^{H}(2_{\fS_3})$ and is a trivial source $kH$-module with vertex $Q$. In addition, the number of isomorphism classes of trivial source $kH$-modules with vertex $P$ is given by the number of $2'$-conjugacy classes of  $N_H(P)/P=P/P\cong \{1\}$, i.e. one. (See e.g. \cite[Theorem~4.6(c)]{LasSurvey}.)
    \par
    Now, as $\bc$ is splendidly Morita equivalent to $\bB_0(kH)$, it follows that the simple $\bc$-modules $1$ and $2$ are trivial source $\bc$-modules with vertex $P$ and $Q$ respectively. Also, since $\bB_0(kH)$  contains, up to isomorphism, a unique trivial source module with vertex $P$, so does $\bc$, as claimed.
\end{enumerate}    
\end{proof}

\begin{lem}\label{lem:tsmInA}
Let $f$ denote the Green correspondence with respect to $(G,L,P)$ from $G$ to $L=N_G(Q)$. Then, the following assertions hold.
\begin{enumerate}[label={\rm(\alph*)}]
    \item The simple $\bA$-module $15$ is a trivial source module with vertex $P$ and $\chi^{}_{\widehat{15}}=\chi^{}_{19}$.
    \item The simple $\bA$-module $6$ is
    a trivial source module with vertex $Q$ and $\chi^{}_{\widehat{6}}=\chi^{}_{17}$. 
    \item We have $f(15)=1$ and $f(6)=2$. 
\end{enumerate}
\end{lem}

\begin{proof}
We prove the three claims together. First, we can read from the ATLAS \cite[p.10]{ATLAS} that $G$ has a maximal subgroup $F\cong C_3\times \fS_5$. Now, letting $k_{\zeta_1}$ denote the simple $kC_3$-module corresponding to the non-trivial character $\zeta_1\in\Irr(C_3)$, the $2$-permutation $kG$-module $X:=\Ind_F^{G}(k_{\zeta_1}\otimes_kk_{\fS_5})$ affords the character 
    \[
        \chi^{}_{\widehat{X}}=\Ind_F^{G}(\zeta_1\otimes 1^{}_{\fS_5})=\chi^{}_{17}+\chi^{}_{19}\,
    \] 
(where the decomposition into irreducibles is easily obtained using GAP 4 \cite{GAP}). Hence, we read from the $2$-decomposition matrix of $\bA$ given in Notation~\ref{nota:Section_3A7}\ref{nota:3A72-decompMatrixA} that $X$ is an $\bA$-module and  has two composition factors, each with multiplicity one,  namely $6$ and $15$.
Moreover, \cite[II Theorem 12.4]{LandrockBook} yields
$\dim_k(\End_{kG}(X))=\langle \chi^{}_{\widehat{X}} ,\chi^{}_{\widehat{X}}\rangle_G=2$, and so $X$ cannot be uniserial of length $2$, proving that $X=15\oplus 6$.
\par
It follows that both the simple $\bA$-modules $15$ and $6$ are trivial source $kG$-modules with $\chi^{}_{\widehat{6}}=\chi^{}_{17}$ and $\chi^{}_{\widehat{15}}=\chi^{}_{19}$.  Clearly, $P$ is a vertex of $15$ since its $k$-dimension is odd. Moreover, $\dim_k(6)=2\cdot 3$ implies that $6$ has vertices of order $4$ (see \cite[Lemma~2.2]{Okuyama91}), and as $\chi^{}_{\widehat{6}}(2A)=2$ and $\chi^{}_{\widehat{6}}(4A)=0$ (see Notation~\ref{nota:Section_3A7}\ref{nota:3A72-decompMatrixA}), by~\cite[II Lemma 12.6(ii)+(iii)]{LandrockBook}, 
the vertices of $6$ cannot be cyclic of order $4$, so they are Klein-four subgroups.  Then, to prove that $Q$ is a vertex of $6$, we prove that $Q_2$ cannot be a vertex of $6$. Indeed,  restricting $\chi^{}_{\widehat{6}}$ to $N_G(Q_2)$  using GAP yields
    \[
            \Res^{G}_{N_G(Q_2)}(\chi^{}_{\widehat{6}})=\eta_1+\eta_2
    \]
with $\eta_1,\eta_2\in\Irr(N_G(Q_2))$ both of degree $3$ and taking values $1$ and $-1$ on the unique conjugacy class of elements of order $4$. It follows that $\Res^{G}_{N_G(Q_2)}(6)$ is indecomposable (indeed, if it decomposed, one of the summands would afford a character with a negative value at some $2$-element, which is impossible for a $2$-permutation module by ~\cite[II Lemma 12.6(ii)]{LandrockBook}). This means that $\Res^{G}_{N_G(Q_2)}(6)$ is not simple, because it has at least two composition factors. Consequently, $Q_2$ cannot be a vertex of $6$ since then its Green correspondent in $N_G(Q_2)$ would have to be simple by Okuyama's Theorem \cite[Theorem~2.3]{Okuyama91}. 
\par
Next, $f(15)$ and $f(6)$ lie in $\bc$ as the Green correspondence and the Brauer correspondence commute. Since $\bc$ contains a unique trivial source module with vertex $P$ by Lemma~\ref{lem:SplendidlyEquivBrauerCorrespondents}(c), clearly $f(15)=1$. Finally, $f(6)$ is also simple by Okuyama's Theorem \cite[Theorem~2.3]{Okuyama91}. Thus, the only possibility is $f(6)=2$.
\end{proof}

\begin{lem}\label{lem:simple-tensor-M}%
Let $M$ be the $(\bc,\bA)$-bimodule defined in Notation~\ref{nota:Section_3A7}{\rm\ref{nota:3A7_candidate_bimodule}}. Then, the following assertions hold:
\begin{enumerate}[label={\rm(\alph*)}]
    \item $M\otimes_{\bA}15\cong 1$ and $M\otimes_{\bA}6\cong 2$ as $\bc$-modules;
    \item $M$ and its dual realise a splendid Morita equivalence between $\bA$ and $\bc$.
\end{enumerate}
\end{lem}

\begin{proof}
By Proposition~\ref{prop:gluing-3A7}, we already know that $M$ realises a stable equivalence of Morita type between $\bA$ and $\bc$. Therefore, we may use \cite[Theorem 2.1(ii) and (iii)]{Lin96StableEqs}. 
\begin{enumerate}[label={\rm(\alph*)}]
    \item First, by \cite[Theorem 2.1(ii)]{Lin96StableEqs}, the $\bc$-modules $M\otimes_{\bA}15$ and $M\otimes_{\bA}6$ are indecomposable and non-projective, since $15$ and $6$ are non-projective simple $\bA$-modules. 
    Moreover, by Lemma~\ref{lem:tsmInA}, the simple modules $15$ and $6$ are trivial source modules with vertex $P$ and $Q$ respectively. Thus, we know from \cite[Lemma A.3(ii) and (iii)]{KMN11} that $M\otimes_{\bA}15$ is a trivial source module with vertex $P$ and 
    $M\otimes_{\bA}6$ is a trivial source module with vertex~$Q$.
    Now, $M=\mathfrak{f}(\bA)|\Res^{G\times G}_{L\times G}(\mathfrak{\bA})$
    implies that 
    \[
    M\otimes_{\bA}S\,|\,\Res^{G\times G}_{L\times G}(\mathfrak{\bA})\otimes_{\bA}S=\Res^G_L(S)
    \]
    for each $S\in \Irr_k(\bA)=\{15,6\}$. 
    By the properties of the Green correspondence, we know that 
    \[
        \Res^G_L(15)=f(15)\oplus X\qquad\text{and}\qquad\Res^G_L(6)=f(6)\oplus Y\,,
    \] 
    where $X$ (respectively $Y$) is a direct sum of indecomposable $kL$-modules with a vertex which is not $L$-conjugate to $P$ (respectively not $L$-conjugate to $Q$). Thus, it follows from Lemma~\ref{lem:tsmInA}(c) that $M\otimes_{\bA}15=f(15)=1$ and $M\otimes_{\bA}6=f(6)=2$. 

    \item By (a), for any simple $\bA$-module $S$, the $\bc$-module $M\otimes_{\bA}S$ is simple.  Therefore, $M$ realises a Morita equivalence between $\bA$ and $\bc$ by \cite[Theorem~2.1(iii)]{Lin96StableEqs}.
    Since by definition $M$ is a  $2$-permutation bimodule, this equivalence is splendid, as required.
\end{enumerate}
\end{proof}

\begin{thm}\label{thm:FaithfulBlocks3A7}%
    The two faithful blocks $\bA,\bA^\ast\in\Bl^{\mathrm{max}}(k[3.\fA_7])$ are splendidly Morita equivalent to the principal block $\bB_0(k\fS_4)=k\fS_4$ of the symmetric group $\fS_4$.
\end{thm}

\begin{proof}
By Lemma~\ref{lem:simple-tensor-M}(b) the block $\bA$ is splendidly Morita equivalent to its Brauer correspondent~$\bc$ in $L=N_G(Q)$. As the Brauer correspondence and the Green correspondence commute with duality, analogous arguments show that $\bA^\ast$ is splendidly Morita equivalent to its Brauer correspondent~$\bc^\ast$. 
Furthermore, by Lemma~\ref{lem:SplendidlyEquivBrauerCorrespondents}, both the blocks $\bc$ and $\bc^\ast$ are splendidly Morita equivalent to $\bB_0(kH)=kH$, where $H\cong \fS_4$. The claim follows. 
\end{proof}

\bigskip

\section{The faithful blocks of maximal defect of \texorpdfstring{$3.\fA_6$}{3.A6} }\label{sec:3A6}%

The group $\PSL_2(9)$ is isomorphic to~$\fA_6$, and its triple cover $3.\fA_6$ possesses two faithful $2$-blocks of maximal defect, denoted Block~4 and Block~5 in \cite[A6mod2]{ModularAtlas}, which are contragredient to each other. 
Below, we prove that these blocks are splendidly Morita equivalent to the principal block of $\PSL_2(7)$. Again, it is clear from \cite{ErdmannTameBook} and \cite[Theorem~2.1]{MacGregor22} that these blocks are Morita equivalent, because they have the same decomposition matrix, but Erdmann's approach does not provide us with a bimodule realising the equivalence. We construct such a bimodule.\par
We note, however, that contrary to the previous case, the Brauer correspondence cannot be used, as it is easily checked using GAP4~\cite{GAP} and \cite[Theorem~2.1]{MacGregor22} that none of the subgroups $L$ of $3.\fA_6$ with $|L|_2=8$ has blocks with the same decomposition matrix. Moreover, not all simple modules involved are trivial source modules. 

\vspace{2mm}
\subsection{Context and notation}
As in Section~\ref{sec:3A7}, we first introduce the notation necessary to solve the problem carefully. 

\begin{nota}\label{nota:Section_3A6}%
Throughout this subsection, we use the following notation. First, 
\begin{enumerate}[label=(\arabic*)]
    \item $G:=3.\fA_7$ as in Section~\ref{sec:3A7}; 
    \item $H:=\PSL_2(7)\cong\GL_3(2)$;
    \item $F:=3.\fA_6$ denotes the triple cover of $\fA_6$;
    \item $Z:=Z(F)=Z(G)=:\langle z\rangle\cong C_3$\,.
\end{enumerate}
We may, and do, assume that $F$ and $H$ are subgroups of $G$. 
\begin{enumerate}[label=(\arabic*), start=5]  
    \item We set $J:= H\cap F\leq G$ and note that $J\cong \fS_4$.
    \item We fix a Sylow $2$-subgroup $P\in \Syl_2(G)$ such that $P\in\Syl_2(H)$ and
    $P\in\Syl_2(F)$\,. 
    We have $P\cong D_8$, $N_F(P)=Z\times P\cong C_3\times D_8$ and $N_H(P)=P\cong D_8$\,.
    \item We let $Q,Q_2\leq P$ denote the two   Klein-four subgroups of $P$ as in Notation~\ref{nota:Section_3A7}.\\
    These are not $F$-conjugate. 
    Notice, however, that in this case $N_F(Q)\cong C_3\times \fS_4$ and also $N_F(Q_2)\cong C_3\times \fS_4$ (where the $C_3$ factor is $Z(F)$) but $N_F(Q)\neq N_F(Q_2)$.   
\end{enumerate}
Thus, we have the following chains of subgroups:
 $$Q< P < N_F(P)=N_{G}(P)< N_F(Q)=N_{G}(Q)< F< G$$
 and 
 \[
 Q< P=N_H(P)< N_H(Q)< H< G\,.
 \]
Next, we consider the following blocks: 
\begin{enumerate}[label=(\arabic*), start=8] 
    \item\label{nota:3A6_B0(kH)} $\bB:=\bB_0(kH)$, the principal block of $kH$;
    \item $\bC:=$ Block 4 of $kF$ according to the notation of \cite[A6mod2]{ModularAtlas},\\
    $\bC^{\ast}:=$ Block 5 of $kF$ according to the notation of \cite[A6mod2]{ModularAtlas}. 
\end{enumerate} 
Without loss of generality we may assume that $\bC$ and $\bC^\ast$ are the Brauer correspondents of the blocks $\bA$ and $\bA^\ast$ of $k[3.\fA_7]$, respectively. 
Moreover, as in Section~\ref{sec:3A7}, we write:
\begin{enumerate}[label=(\arabic*), start=10]
    \item\label{nota:blocksOfkZ_3A6}  $\Irr(Z):=\{\zeta_0:=1_{C_3},\zeta_1,\zeta_2\}$, $e_i:=e_{\zeta_i}\in kZ$ for the block idempotent associated to $\zeta_i$ for each $0\leq i\leq 2$, hence  $kZ=kZe_0\oplus kZe_1\oplus kZe_2$ is the block decomposition of $kZ$. 
 \end{enumerate} 
We have that $P$ is a defect group of $\bC$, $\bC^\ast$ and $\bB$, whereas the blocks of $kZ$ are of defect zero.\\ 
The simple $\bC$-modules and the simple $\bB$-modules are as follows.
\begin{enumerate}[label=(\arabic*),start=11]
    \item\label{nota:3A6-decompMatrixB} By~\cite[L3(2)mod2]{ModularAtlas} the principal block $\bB$ of $k\PSL_2(7)$ has three simple modules: the trivial module $k_H$ and two simple modules of dimension $3$.   Following~\cite[L3(2)mod2]{ModularAtlas} we set $\Irr_k(\bB):=\{k_H, \tilde{3}_1,\tilde{3}_2\}$\,.  
    The $2$-decomposition matrix of $\bB$ and the values at the conjugacy classes of $2$-elements are: 
        \[
        \begin{tabular}{c:c|ccc}{\footnotesize{\cite[p.3]{ATLAS}}} &  \cite{ModularAtlas} & $\tilde{3}_1$ & $\tilde{3}_2$ & $k_H$\\
            \hline  
            $\chi^{}_{2}$ & $\chi^{}_{3_1}$  & $1$ & $0$ & $0$\\
            $\chi^{}_{3}$ & $\chi^{}_{3_2}$ & $0$ & $1$ & $0$\\
            $\chi^{}_{4}$ & $\chi^{}_{6_1}$ & $1$ & $1$ & $0$\\
            $\chi^{}_{1}$ & $\chi^{}_{1_1}$ & $0$ & $0$ & $1$\\
            $\chi^{}_{5}$ & $\chi^{}_{7_1}$ & $1$ & $1$ & $1$\\
        \end{tabular}
        \qquad \qquad
        \begin{tabular}{c|ccc} &  $1A$ & $2A$ & $4A$\\
            \hline  
            $\chi^{}_{2}$ & $3$  & $-1$ & $1$\\
            $\chi^{}_{3}$ & $3$ & $-1$ & $1$\\
            $\chi^{}_{4}$ & $6$ & $2$ & $0$\\
            $\chi^{}_{1}$ & $1$ & $1$ & $1$\\
            $\chi^{}_{5}$ & $7$ & $-1$ & $-1$\\
        \end{tabular}
        \]
    \item\label{nota:3A6-decompMatrixA} By~\cite[A6mod2]{ModularAtlas} the blocks $\bC$ and $\bC^\ast$ both have three simple modules: two of dimension 3 and one of dimension 9. Following~\cite[A6mod2]{ModularAtlas} we set $\Irr_k(\bC):=\{3_1,3_3, 9_1\}$  and $\Irr_k(\bC^\ast):=\{3_1^\ast, 3_3^\ast, 9_1^\ast\}$. The $2$-decomposition matrix of $\bC$ and the values at the conjugacy classes of $2$-elements are: 
        \[
        \begin{tabular}{c:c|ccc}{\footnotesize{\cite[p.4]{ATLAS}}} &  \cite{ModularAtlas} & $3_1$ & $3_3$ & $9_1$\\
            \hline  
            $\chi^{}_{14}$ & $\chi^{}_{3_1}$  & $1$ & $0$ & $0$\\
            $\chi^{}_{15}$ & $\chi^{}_{3_3}$ & $0$ & $1$ & $0$\\
            $\chi^{}_{16}$ & $\chi^{}_{6_1}$ & $1$ & $1$ & $0$\\
            $\chi^{}_{17}$ & $\chi^{}_{9_2}$ & $0$ & $0$ & $1$\\
            $\chi^{}_{18}$ & $\chi^{}_{15_1}$ & $1$ & $1$ & $1$\\
        \end{tabular}
        \qquad \qquad
        \begin{tabular}{c|ccc} &  $1A$ & $2A$ & $4A$\\
            \hline  
            $\chi^{}_{14}$ & $3$  & $-1$ & $1$\\
            $\chi^{}_{15}$ & $3$ & $-1$ & $1$\\
            $\chi^{}_{16}$ & $6$ & $2$ & $0$\\
            $\chi^{}_{17}$ & $9$ & $1$ & $1$\\
            $\chi^{}_{18}$ & $15$ & $-1$ & $-1$\\
        \end{tabular}
        \]        
\end{enumerate} 
\end{nota}

\vspace{2mm}
\subsection{A stable equivalence of Morita type}
As in Subsection~\ref{ssec:stEqMorTp} we first prove a stable equivalence of Morita type between $\bC$ and $\bB$ using Linckelmann's glueing method~\cite[Theorem~3.1]{Lin01}. For this purpose, we collect information about the fusion systems of $\bC$ and $\bB$ and about the structure of the centralisers of the $2$-subgroups, which can be found at \cite[Groups/Abstract/1080.260 and Groups/Abstract/168.42]{lmfdb} or using GAP4~\cite{GAP}.

\begin{nota}\label{nota:Section_3A6_part_II}\leavevmode 
\begin{enumerate}[label=(\arabic*),start=13]
    \item\label{nota:local-3A6} 
    We have $C_H(P)=Z(P)\cong C_2$. It follows that $kC_H(P)\cong kC_2$ has just one block idempotent, hence $\Br^H_P(e_{\bB})=1_{kC_H(P)}=:e_P$.\\
    We have $C_F(P)=Z\times Z(P)\cong C_3\times C_2$ and, under the identification
        \[
            k[Z\times Z(P)]\cong kZ\otimes_k kZ(P)\,,
        \]
    we have
        \(
            kC_F(P)\Br^F_P(e_\bC)=k[Z\times Z(P)](e_1\otimes 1_{kZ(P)})\,.
\)
(Up to permuting the labelling of $e_1$ and $e_2$ if necessary.)
We identify $e_1$ with $e_1\otimes 1_{kZ(P)}$ and we set 
\[
f_P:=\Br^F_P(e_\bC)=e_1\,.
\]
In particular, $(P,f_P)$ and $(P,e_P)$ are maximal Brauer pairs of $\bC$ and $\bB$, respectively.

\item For any $R\le P$, let $e_R$ (resp. $f_R$) denote the unique block idempotent of $kC_H(R)$ (resp. $kC_F(R)$) such that 
\[(R,e_R)\le_H(P,e_P)\qquad\text{and}\qquad (R,f_R)\le_F (P,f_P)\,.
\] 
The fusion systems of $\bB$ and $\bC$ with respect to the maximal Brauer pairs $(P,e_P)$ and $(P,f_P)$ are isomorphic, 
as both blocks have the same number of isomorphism classes of simple modules. 
(This follows from~\cite[Part IV, Section 5.7]{AschKessOli11} together
with~\cite{BrauerDihedral}; see also \cite[Theorem 8.1(2)]{Sambale14}.)
We identify them and denote the resulting fusion system on $P$ by~$\calF$. In this case, $P$ has one $\calF$-conjugacy class of involutions, namely the class of the central involution. 
\item \label{rem:local-3A6}
     If $R\le P$ is one of $P$, $Q$, $Q_2$, $Z(P)$ or the cyclic group of order $4$, then one has
    \[
    C_H(R)=C_P(R)\qquad\text{and}\qquad C_F(R)=Z\times C_P(R)\,.
    \]
    Since $C_P(R)$ is a $2$-group, under the natural identification $kC_F(R)\cong kZ\otimes_k kC_P(R)$, we have $f_R=e_1\otimes 1_{kC_P(R)}$. As before, identifying $e_1$ with $e_1\otimes 1_{kC_P(R)}$, we simply write $f_R=e_1$.
\end{enumerate}
\end{nota}    


\begin{rem}\label{rem:diagonal-normalizers-3A6} Let $R\le P$. Then
\[
C_H(R)\times C_F(R)=C_{H\times F}(\Delta R)
\leq 
N_{H\times F}(\Delta R)
\leq N_H(R)\times N_F(R)\,.
\]
The latter inclusion is an equality if and only if both $N_H(R)$ and $N_F(R)$ centralise $R$. In particular, if $R\cong C_2$, we have 
\[
C_{H\times F}(\Delta R)=C_H(R)\times C_F(R)=N_H(R)\times N_F(R)=N_{H\times F}(\Delta R)\,.
\]
For any $R\leq P$ and any $k[H\times F]$-module $Y$, the Brauer
construction $Y[\Delta R]$ is naturally a
$kN_{H\times F}(\Delta R)$-module. By restriction to $C_{H\times F}(\Delta R)=C_H(R)\times C_F(R)$, we may therefore regard $Y[\Delta R]$ as a
$(kC_H(R),kC_F(R))$-bimodule, with bimodule structure given by
\[
h\cdot y\cdot f:=(h,f^{-1})y
\]
for every $h\in C_H(R)$, $f\in C_F(R)$ and $y\in Y[\Delta R]$.
\end{rem}

\noindent Next we define a series of bimodules, which will help us prove the required Morita equivalence between $\bB$ and $\bC$.

\begin{nota}\leavevmode
\begin{enumerate}[label=(\arabic*),start=16]
   \item Set  $\widetilde{X}:=\Ind_{\Delta P}^{H\times F}(k)$, and consider the $(\bB,\bC)$-bimodule $X:=e_{\bB}\widetilde{X}e_{\bC}$ 
   \item  Set $\widetilde{W}:=\Ind^{H\times F}_{\Delta J}(k)$ and  consider the $(\bB,\bC)$-bimodule $W:=e_\bB \widetilde{W} e_{\bC}$.
\end{enumerate}    
\end{nota}

\begin{lem}\label{lem:Brauer-fixed-points-W}
Let $R$ be one of $P$, $Q$, $Q_2$ or the cyclic group of order $4$. Then
\[
\widetilde{W}[\Delta R]\cong k[N_{H\times F}(\Delta R)/\Delta N_J(R)]
\]
as $kN_{H\times F}(\Delta R)$-modules, where the left action on the permutation module on the right-hand side is induced by left multiplication.
\end{lem}
\begin{proof}
    First note that $\widetilde{W}$ is isomorphic to the permutation $k[H\times F]$-module $k\Omega$, where $\Omega:=(H\times F)/\Delta J$. Hence, by \cite[Proposition~5.8.1]{LinckelmannVol1}, $\widetilde{W}[\Delta R]\cong k[\Omega^{\Delta R}]$ as $kN_{H\times F}(\Delta R)$-modules, where $\Omega^{\Delta R}$ denotes the set of $\Delta R$-fixed points in $\Omega$, and we identify each element of $\Omega^{\Delta R}$ with its image under the Brauer map.

    A coset $(h,f)\Delta J$ is $\Delta R$-fixed if and only if 
    \[
    h^{-1}rh=f^{-1}rf\in J\qquad \text{for every } r\in R\,.
    \]
    For each of the subgroups $R$ under consideration, every subgroup of
    $J$ which is $F$-conjugate to $R$ is $J$-conjugate to $R$. Hence, in this case, $h^{-1}Rh=f^{-1}Rf$ is $J$-conjugate to $R$. Choose $x\in J$ such that $x^{-1}h^{-1}Rhx=R$. Then
    \[
    (h,f)\Delta J=(hx,fx)\Delta J\qquad \text{and} \qquad (hx,fx)\in N_{H\times F}(\Delta R)\,.
    \]
    Conversely, every coset represented by an element of $N_{H\times F}(\Delta R)$ is $\Delta R$-fixed. Thus $\Omega^{\Delta R}$ is the $N_{H\times F}(\Delta R)$-orbit of $\Delta J$. Since the stabilizer of the coset $\Delta J$ in $N_{H\times F}(\Delta R)$ is
    \[
    N_{H\times F}(\Delta R)\cap \Delta J=\Delta N_J(R)\,,
    \]
    we obtain 
    \[
    \Omega^{\Delta R}\cong N_{H\times F}(\Delta R)/\Delta N_J(R)
    \]
    as $N_{H\times F}(\Delta R)$-sets. Linearizing this isomorphism proves the result.
\end{proof}

\begin{lem}\label{lem:candidate_3A6}
The following assertions hold. 
\begin{enumerate}[label=\rm(\alph*)]
    \item We have $N_{H\times F}(\Delta P)=\Delta P\cdot({1}\times (Z\times Z(P)))$ and  the map 
    \[
    N_{H\times F}(\Delta P)/\Delta P\longrightarrow Z\times Z(P), \qquad (a,b)\Delta P\mapsto ab^{-1}
    \]
    is an isomorphism.
    \item  We have $\widetilde{W}[\Delta P]\cong k[Z\times Z(P)]$ and $W[\Delta P]\cong k[Z\times Z(P)]e_1$ as $kN_{H\times F}(\Delta P)$-modules. Here the action of $N_{H\times F}(\Delta P)$ on $Z\times Z(P)$ is the transport, via the isomorphism in part~{\rm (a)}, of the left multiplication action on $N_{H\times F}(\Delta P)/\Delta P$. More explicitly, 
    \[
    (u,v)\cdot x=uxv^{-1}
    \]
    for $(u,v)\in N_{H\times F}(\Delta P)$ and $x\in Z\times Z(P)$.
    In particular, $W[\Delta P]$ is indecomposable.
    \item The bimodule $X$,  seen as a $k[H\times F]$-module, has a unique indecomposable summand with vertex $\Delta P$, say~$M$. In particular, $M$ is, up to isomorphism, the unique $(kH,kF)$-bimodule with vertex $\Delta P$ that can realise a splendid Morita equivalence between $\bB$ and $\bC$. 
    \item The $k[H\times F]$-module $W$ has a unique direct summand with vertex $\Delta P$, and this summand is isomorphic to the $k[H\times F]$-module $M$ in Part~{\rm (c)}.
\end{enumerate}
\end{lem}

\begin{proof}
   \begin{enumerate}[label=\rm(\alph*)]
    \item Since $N_F(P)=Z\times P$ and $N_H(P)=P$, intersecting $N_{G\times G}(\Delta P)$ with $H\times F$, \cite[Proposition 2.4(a)]{BoltjePerepelitsky} gives
    \[
    N_{H\times F}(\Delta P)
    =\Delta P\cdot(\{1\}\times C_F(P))
    =\Delta P\cdot(\{1\}\times(Z\times Z(P)))\,.
    \] 
The second claim is immediate after noting that $\Delta P\cap (\{1\}\times(Z\times Z(P)))=1$.
\item Since $N_J(P)=P$, by Lemma~\ref{lem:Brauer-fixed-points-W}, $\widetilde{W}[\Delta P]\cong k[N_{H\times F}(\Delta P)/\Delta P]$. Transporting the left multiplication action through the isomorphism in part~{\rm (a)} gives
\[
\widetilde{W}[\Delta P]\cong k[Z\times Z(P)]
\]
as $kN_{H\times F}(\Delta P)$-modules. Indeed, if $(u,v)\in N_{H\times F}(\Delta P)$ and $x\in Z\times Z(P)$, then the
transported action is given by $(u,v)\cdot x=uxv^{-1}$.

By \cite[Lemma~3.7]{BoltjePerepelitsky},
\[
W[\Delta P]\cong\Br_P(e_\bB)\widetilde{W}[\Delta P]\Br_P(e_\bC)\cong k[Z\times Z(P)]e_1.
\]
Finally, upon restriction to $Z(P)\times 1\leq C_H(P)\times C_F(P)$, the latter module is isomorphic to the regular $kZ(P)$-module which is indecomposable. Hence $W[\Delta P]$ is indecomposable.

\item Since $\widetilde{X}$ is the permutation $k[H\times F]$-module $k[(H\times F)/\Delta P]$, we have  
    \[
    \widetilde{X}[\Delta P]\cong \Ind_{\Delta P}^{N_{H\times F}(\Delta P)}(k)\cong k[N_{H\times F}(\Delta P)/\Delta P]
    \]
    as $k[N_{H\times F}(\Delta P)]$-modules; see e.g. \cite[Corollary 5.8.4]{LinckelmannVol1}.
Therefore, by part~{\rm (a)} and \cite[Lemma 3.7]{BoltjePerepelitsky}, we have 
\[
X[\Delta P]\cong \Br_P(e_\bB)\widetilde{X}[\Delta P]\Br_P(e_\bC)\cong k[Z\times Z(P)]e_1\,.
\]
By part~{\rm (b)}, the latter module is indecomposable and the claim follows from the Brou\'{e} correspondence; see~\cite[Theorem (3.2) and (3.4)]{BrouePperm85}. 
\item By part~{\rm (b)}, $W$ has a unique
indecomposable summand, say $M'$, with vertex $\Delta P$. Hence
\[
M'\mid \Ind_{\Delta P}^{H\times F}(k)=\widetilde X.
\]
Moreover, since $M'$ is a direct summand of
$W=e_\bB\widetilde W e_\bC$, the block idempotents $e_\bB$ and $e_\bC$
act as the identity on $M'$. Thus $M'$ is in fact a direct summand of $X=e_\bB\widetilde X e_\bC$. By part~{\rm (c)}, $M'\cong M$, as required.
    \end{enumerate}
\end{proof}

\begin{nota}\label{nota:BimodulesSection_3A6}\leavevmode
\begin{enumerate}[label=(\arabic*),start=18]
    \item\label{nota:DefBimodule-3A6} 
    Let $M$ denote the unique indecomposable direct summand of $X$ with vertex $\Delta P$ given by  Lemma~\ref{lem:candidate_3A6}(c), which we may see as a $(\bB,\bC)$-bimodule.  
\end{enumerate}    
\end{nota}

\begin{lem}\label{lem:non-vanishing-local-M}
    For every subgroup $R\le P$, one has $e_RM[\Delta R]f_R\neq 0$.
\end{lem}
\begin{proof}
    By Lemma~\ref{lem:candidate_3A6}(b),
    \[
    e_PW[\Delta P]f_P=W[\Delta P]\cong k[Z\times Z(P)]e_1\neq 0\,.
    \]
    By Lemma~\ref{lem:candidate_3A6}(d), $M$ is, up to isomorphism, the unique direct summand of $W$ with vertex $\Delta P$. It follows that $e_PM[\Delta P]f_P\neq 0$. Thus, by \cite[Proposition~5.3(c)]{BoltjePerepelitsky}, $(\Delta P, e_P\otimes f_P^*)$ is a maximal $M$-Brauer pair. By \cite[Remark~10.2(c)]{BoltjePerepelitsky}, we have $(\Delta R, e_R\otimes f_R^*)\le (\Delta P,e_P\otimes f_P^*)$. Hence, by \cite[Proposition~5.3(b)]{BoltjePerepelitsky}, the pair $(\Delta R, e_R\otimes f_R^*)$ is an $M$-Brauer pair, as well. Thus, $e_RM[\Delta R]f_R$ is non-zero.
\end{proof}

\begin{lem}\label{lem:Q_2-local-equiv}
    Let $R=Q_2$. Then $e_RM[\Delta R]f_R\cong k[Z\times R]e_1$ as $(kC_H(R)e_R, kC_F(R)f_R)$-bimodules. In particular, $e_RM[\Delta R]f_R$ induces a splendid Morita equivalence between $kC_H(R)e_R$ and $kC_F(R)f_R$.
\end{lem}
\begin{proof}
    First note that $N_H(R)\cong \fS_4$ and $N_F(R)\cong Z\times \fS_4$. Choose an involution $t\in P\setminus R$, so that $P=R\rtimes \langle t\rangle$. Choose complements $S_H$ to $R$ in $N_H(R)$ and $S_F$ to $Z\times R$ in $N_F(R)$, both containing $t$. Note that both $S_H$ and $S_F$ are isomorphic to $\fS_3$. Let $\theta:S_H\longrightarrow S_F$ be an isomorphism such that $s$ and $\theta(s)$ induce the same automorphism of $R$ for every $s\in S_H$. In particular, $\theta(t)=t$. Then
    \[
    N:=N_{H\times F}(\Delta R)=\{(as,zb\theta(s))\in N_H(R)\times N_F(R)\mid a,b\in R, z\in Z, s\in S_H\}\,. 
    \]
    Define 
    \[
    \Phi:N\longrightarrow Z\times N_H(R), \qquad (as, zb\theta(s))\mapsto (z^{-1}, ab^{-1}s)\,.
    \]
    A direct calculation shows that $\Phi$ is a surjective homomorphism with $\ker(\Phi)=\Delta R$ and $\Phi(\Delta P)=1\times\langle t\rangle$. Consider $(Z\times N_H(R)/(1\times\langle t\rangle))$ as an $N$-set via $\Phi$, that is, for $n\in N$ and $x\in Z\times N_H(R)$, we set
    \[
    n\cdot x(1\times \langle t\rangle)=\Phi(n)x(1\times \langle t\rangle)\,.
    \]
    Since $\ker(\Phi)=\Delta R\le \Delta P$, we have $\Phi^{-1}(1\times \langle t\rangle)=\Delta P$. Thus, $\Phi$ induces an isomorphism of $N$-sets
    \[
    \overline{\Phi}: N/\Delta P\longrightarrow (Z\times N_H(R))/(1\times \langle t\rangle), \qquad n\Delta P\mapsto \Phi(n)(1\times \langle t\rangle)\,.
    \]
    Using the identification $(Z\times N_H(R))/(1\times \langle t\rangle)\cong Z\times (N_H(R)/\langle t\rangle)$, we obtain an isomorphism of $N$-sets
    \[
    N/\Delta P\cong Z\times (N_H(R)/\langle t\rangle)\,.
    \]
    Since $N_J(R)=P$, together with Lemma~\ref{lem:Brauer-fixed-points-W}, this implies that we have isomorphisms
    \[
    \widetilde{W}[\Delta R]\cong k[N/\Delta P]\cong k[Z\times (N_H(R)/\langle t\rangle)]
    \]
    of $kN$-modules. Therefore,
    \[
    \widetilde{W}[\Delta R]\cong kZ\otimes \Ind_{\langle t\rangle}^{N_H(R)}(k)
    \]
    as $kN$-modules, where the $N$-action on the latter module is induced by $\Phi$. More precisely, for every $n=(as,zb\theta(s))\in N$, $z_0\in Z$ and $x\in N_H(R)$,
    \[
    n\cdot (z_0\otimes (x\otimes 1))=z^{-1}z_0\otimes (ab^{-1}sx\otimes 1)\,.
    \]
    We now decompose the $kN_H(R)$-module $\Ind^{N_H(R)}_{\langle t\rangle} k$. Since $\langle t\rangle\le S_H\le N_H(R)$, we get 
    \[
    \Ind^{N_H(R)}_{\langle t\rangle}(k)\cong \Ind^{N_H(R)}_{S_H}(\Ind^{S_H}_{\langle t\rangle}(k))\,.
    \]
    But $\Ind^{S_H}_{\langle t\rangle}(k)\cong k\oplus E$ as $kS_H$-modules, where $E$ is the two-dimensional simple projective $kS_H$-module. Hence $\Ind^{N_H(R)}_{\langle t\rangle}(k)\cong U_0\oplus V_0$, where $U_0=\Ind^{N_H(R)}_{S_H}(k)$ and $V_0=\Ind^{N_H(R)}_{S_H}(E)$. The module $U_0$ is indecomposable with vertex $\langle t\rangle$. Indeed, $\Res^{N_H(R)}_R(U_0)\cong kR$ is indecomposable, and thus $U_0$ is also indecomposable. Also, $k$ is a direct summand of $\Ind^{S_H}_{\langle t\rangle}(k)$, so the transitivity of induction yields $U_0\mid \Ind^{N_H(R)}_{\langle t\rangle}(k)$. Thus, $U_0$ is relatively $\langle t\rangle$-projective. Since $\dim_k(U_0)=4$ is not divisible by $|N_H(R)|_2=8$, $U_0$ is not projective. Thus, $\langle t\rangle$ is a vertex of $U_0$. Since induction preserves projectivity, $V_0$ is projective. Since $kZe_1$ is a one-dimensional projective as a $kZ$-module, the $k[Z\times N_H(R)]$-module $kZe_1\otimes_k U_0$ is indecomposable with vertex $1\times \langle t\rangle$. In addition, $kZe_1\otimes_k V_0$ is projective. 

    Recall that $e_R=1$ and $f_R=e_1$. Thus,
    \[
    e_R\widetilde{W}[\Delta R]f_R\cong kZe_1\otimes_k \Ind^{N_H(R)}_{\langle t\rangle}(k)\cong U\oplus V
    \]
    where $U=kZe_1\otimes_k U_0$ and $V=kZe_1\otimes V_0$ are viewed as $kN$-modules through $\Phi$. Note that the indecomposable $kN$-module $U$ has vertex $\Phi^{-1}(1\times \langle t\rangle)=\Delta P$ and every indecomposable direct summand of $V$ has vertex $\ker(\Phi)=\Delta R$. 

    Since $N_{\Delta P}(\Delta R)= \Delta P$, $\Delta R$ is fully normalized in $\calF_{\Delta P}(H\times F)$. Since also $M$ is an indecomposable trivial source $k[H\times F]$-module with vertex $\Delta P$, \cite[Lemma~6.4]{KunugiSuzuki} implies that every indecomposable direct summand of $M[\Delta R]$, as a $kN$-module, has a vertex properly containing $\Delta R$. In particular, the same holds for every indecomposable direct summand of $e_RM[\Delta R]f_R$. By Lemma~\ref{lem:candidate_3A6}(d), $M$ is isomorphic to a direct summand of $W$, hence $e_RM[\Delta R]f_R$ is isomorphic to a direct summand of $e_R\widetilde{W}[\Delta R]f_R$. By Lemma~\ref{lem:non-vanishing-local-M}, this summand is non-zero. It follows from the preceding paragraph that $e_RM[\Delta R]f_R\cong U$ as $kN$-modules. 

    Now we restrict the $kN$-module $U=kZe_1\otimes_k U_0$ to $C:=C_{H\times F}(\Delta R)=C_H(R)\times C_F(R)=R\times (Z\times R)$. For $(a,zb)\in R\times (Z\times R)$, we have $\Phi(a,zb)=(z^{-1}, ab^{-1})$. Thus,
    \[
    \Res^N_C(U)\cong kZe_1\otimes_k \Res^{N_H(R)}_R(U_0)\cong kZe_1\otimes kR
    \]
    as $kC$-modules. Here, $kC$ acts on the latter module via $(a,zb)\cdot (u\otimes y)=z^{-1}u\otimes ab^{-1}y$ for every $a,b\in R, z\in Z$, $u\in kZe_1$ and $y\in kR$. The standard $k$-algebra isomorphism
    \[
    \psi:kZe_1\otimes_k kR\longrightarrow k[Z\times R]e_1, \qquad u\otimes x\mapsto ux
    \]
    is in fact an isomorphism of $kC$-modules where $kC$ acts on $k[Z\times R]e_1$ by $(a,zb)\cdot w=aw(zb)^{-1}$ for every $a,b\in R$, $z\in Z$ and $w\in k[Z\times R]e_1$. Indeed, for $a,b\in R$, $z\in Z$, $u\in kZe_1$ and $y\in kR$, we have
    \[
    \psi((a,zb)\cdot (u\otimes y))\!=\psi(z^{-1}u\otimes ab^{-1}y)=z^{-1}uab^{-1}y=a(uy)(zb)^{-1}\!=(a,zb)\cdot (uy)=(a,zb)\cdot\psi(u\otimes y).
    \]
    By Remark~\ref{rem:diagonal-normalizers-3A6}, this $kC$-module structure on $k[Z\times R]e_1$ corresponds to the usual $(kR, k[Z\times R])$-bimodule structure given by left and right multiplication. This proves both assertions.
    \end{proof}

\begin{prop}\label{prop:gluing-3A6}
The $(\bB,\bC)$-bimodule $M$  defined in Notation~\ref{nota:BimodulesSection_3A6}{\rm\ref{nota:DefBimodule-3A6}} and its dual 
induce a stable equivalence of Morita type between the blocks $\bB$ and $\bC$. 
\end{prop}
\begin{proof}
By \cite[Theorem~3.1]{Lin01} and the isomorphism of the fusion systems of $\bB$ and $\bC$, it suffices to show that the bimodule $e_RM[\Delta R]f_R$ induces a splendid Morita equivalence between $kC_H(R)e_R$ and $kC_F(R)f_R$ for representatives $R$ of the non-trivial $\calF$-conjugacy classes of subgroups of $P$. We may take $R\in\{P,Q,Q_2,Z(P),C_4\}$,
where $C_4$ denotes the cyclic subgroup of $P$ of order $4$. By Lemma~\ref{lem:non-vanishing-local-M}, all bimodules $e_RM[\Delta R]f_R$ are non-zero.

\medskip 
\noindent{{\bf Case \rm(i)}.} Suppose that $R=P$. By Lemma~\ref{lem:candidate_3A6}(d),
$M$ is isomorphic to a direct summand of $W$. Hence $e_PM[\Delta P]f_P$ is isomorphic to a non-zero direct summand of $W[\Delta P]$. By Lemma~\ref{lem:candidate_3A6}(b), the latter module is indecomposable, and therefore $e_PM[\Delta P]f_P\cong W[\Delta P] \cong k[Z\times Z(P)]e_1$. By Remark~\ref{rem:diagonal-normalizers-3A6}, this is an isomorphism of
$(kC_H(P)e_P,kC_F(P)f_P)$-bimodules. Thus the required splendid Morita
equivalence follows.

\medskip
\noindent{{\bf Case \rm(ii)}.} Suppose that $R=Q$ or $R\cong C_4$. In both cases, $C_H(R)=R$, $C_F(R)=Z\times R$ and $N_{H\times F}(\Delta R)=\Delta N_J(R)\cdot (\{1\}\times (Z\times R))$. Indeed, if $(a,b)\in N_{H\times F}(\Delta R)$, then $a\in N_J(R)$ and $a^{-1}b\in C_F(R)$. The map
\[
N_{H\times F}(\Delta R)/\Delta N_J(R)\longrightarrow Z\times R,\qquad (a,b)\Delta N_J(R)\mapsto ab^{-1},
\]
is an isomorphism of $N_{H\times F}(\Delta R)$-sets, where the action on the latter set is given by $(u,v)\cdot x=uxv^{-1}$ for every $(u,v)\in N_{H\times F}(\Delta R)$ and $x\in Z\times R$. Thus Lemma~\ref{lem:Brauer-fixed-points-W} gives $\widetilde{W}[\Delta R]\cong k[Z\times R]$ as $k[N_{H\times F}(\Delta R)]$-modules. By \cite[Lemma~3.7]{BoltjePerepelitsky},
\[
W[\Delta R]\cong \Br_R(e_\bB)\widetilde{W}[\Delta R]\Br_R(e_\bC)\cong k[Z\times R]e_1
\]
as $k[N_{H\times F}(\Delta R)]$-modules. The latter module is indecomposable. By Lemma~\ref{lem:candidate_3A6}(d) and Lemma~\ref{lem:non-vanishing-local-M}, $e_RM[\Delta R]f_R$ is isomorphic to a non-zero direct summand of $e_RW[\Delta R]f_R$. Hence $e_RM[\Delta R]f_R\cong k[Z\times R]e_1$. By Remark~\ref{rem:diagonal-normalizers-3A6}, this is an isomorphism of $(kC_H(R)e_R, kC_F(R)f_R)$-bimodules. The latter bimodule induces the desired splendid Morita equivalence.

\medskip
\noindent{{\bf Case \rm(iii).}} The case $R=Q_2$ follows from Lemma~\ref{lem:Q_2-local-equiv}.

\medskip
\noindent{{\bf Case \rm(iv).}} Finally, Suppose that $R=Z(P)$. Since $R\cong C_2$, we have $N_H(R)=C_H(R)=P$, $N_F(R)=C_F(R)=Z\times P$ and
\[
N_{H\times F}(\Delta R)=C_H(R)\times C_F(R)=P\times (Z\times P)\,.
\]
Thus, the $k[N_{H\times F}(\Delta R)]$-module $M[\Delta R]$ is already a $k[C_H(R)\times C_F(R)]$-module. Since $e_R=1_{kP}$ and $f_R=e_1$, the  $(kC_H(R),kC_F(R))$-bimodule $Y:=e_RM[\Delta R]f_R$ is naturally a $(kP,k[Z\times P]e_1)$-bimodule. We first show that $Y$ is indecomposable. So write 
\[
Y=Y_1\oplus\cdots\oplus Y_t\qquad (t\in\mathbb{Z}_{>0})
\]
as a direct sum of indecomposable $k[N_{H\times F}(\Delta R)]$-modules. Since $N_{\Delta P}(\Delta R)=\Delta P$, $\Delta R$ is fully normalized in $\calF_{\Delta P}(H\times F)$. Thus by \cite[Lemma~6.4]{KunugiSuzuki}, every indecomposable direct summand of the $k[N_{H\times F}(\Delta R)]$-module $M[\Delta R]$ has a vertex $S$ satisfying
\begin{equation}\label{eq:vertices-ZP}
\Delta R\lhd S\le_{N_{H\times F}(\Delta R)}\Delta P.
\end{equation}
Since $Y\mid M[\Delta R]$, this conclusion holds in particular for $Y_i$ for each $1\le i\le t$, which can therefore only have vertices of order $4$ or $8$.

Now as $R=Z(P)$ is characteristic in $P$, $N_{H\times F}(\Delta P)\le N_{H\times F}(\Delta R)$. Hence, the transitivity of the Brauer construction \cite[Proposition~5.8.5]{LinckelmannVol1} yields $(M[\Delta R])[\Delta P]\cong M[\Delta P]$ as $k[N_{H\times F}(\Delta P)]$-modules. By the structures of the centralisers described above, we obtain $\Br^{C_H(R)}_P(e_R)=e_P$ and $\Br^{C_F(R)}_P(f_R)=f_P$. Since $e_R$ and $f_R$ are $P$-stable, \cite[Lemma~3.7]{BoltjePerepelitsky} gives
\[
Y[\Delta P]\cong \Br_P^{C_H(R)}(e_R)\left(M[\Delta R]\right)[\Delta P]\Br_P^{C_F(R)}(f_R)\cong e_PM[\Delta P]f_P
\]
as $(kC_H(P)e_P, kC_F(P)f_P)$-bimodules, which is non-zero and indecomposable by Case~{\rm (i)}. On the other hand, the Brauer construction preserves direct sums, so $Y[\Delta P]=Y_1[\Delta P]\oplus\cdots\oplus Y_t[\Delta P]$. It follows that, after relabelling if necessary,
\[
Y_1[\Delta P]\neq 0\qquad\text{and}\qquad Y_i[\Delta P]=0 \qquad \text{for every } 2\le i\le t\,.
\]
As $Y_1$ is a trivial source module, the non-vanishing of $Y_1[\Delta P]$ implies that $\Delta P$ is contained in a vertex of $Y_1$. By \eqref{eq:vertices-ZP}, such a vertex has order at most $|P|$. Hence $\Delta P$ itself is a vertex of $Y_1$. Now fix an index $i$ such that $2\le i\le t$. As $Y_i$ is also a trivial source module, the vanishing of $Y_i[\Delta P]$ shows that $Y_i$ cannot have a vertex which is $N_{H\times F}(\Delta R)$-conjugate to $\Delta P$ and so (\ref{eq:vertices-ZP}) yields that it has a vertex which is $N_{H\times F}(\Delta R)$-conjugate to $\Delta T$, for some subgroup $T\le P$ satisfying $R<T<P$. Therefore, $|T|=4$ and we may assume that $\Delta T$ itself is a vertex of $Y_i$. Again, because $Y_1$ and $Y_i$ are trivial source modules, we have $Y_1[\Delta T]\neq 0$ and $Y_i[\Delta T]\neq 0$, and it follows that $Y[\Delta T]=Y_1[\Delta T]\oplus\cdots\oplus Y_t[\Delta T]$ is decomposable. On the other hand, since $\Delta R\lhd\Delta T$, the transitivity of the Brauer construction gives $Y[\Delta T]\cong e_TM[\Delta T]f_T$ as $(kC_H(T)e_T,kC_F(T)f_T)$-bimodules, where the latter bimodule is non-zero and indecomposable by Cases~{\rm (ii)} and {\rm (iii)}, contradicting the decomposability of $Y[\Delta T]$. Therefore, $t=1$ and $Y$ is indecomposable.

It remains to identify $Y$. The coset $\Delta J\in (H\times F)/\Delta J$ is fixed by $\Delta R$, and its stabilizer in $N_{H\times F}(\Delta R)$ is
\[
N_{H\times F}(\Delta R)\cap\Delta J=\Delta C_J(R)=\Delta P.
\]
Hence by \cite[Proposition~5.8.1]{LinckelmannVol1}, the $N_{H\times F}(\Delta R)$-orbit of $\Delta J$ yields a direct summand $k[N_{H\times F}(\Delta R)/\Delta P]$ of $\widetilde W[\Delta R]$ as $k[N_{H\times F}(\Delta R)]$-modules. Multiplying by $e_R=1_{kP}$ and $f_R=e_1$, we obtain a direct summand $k[N_{H\times F}(\Delta R)/\Delta P]e_1$ of $W[\Delta R]$ as $k[N_{H\times F}(\Delta R)]$-modules.

The map
\[
N_{H\times F}(\Delta R)/\Delta P\longrightarrow Z\times P,\qquad (a,zb)\Delta P\mapsto a(zb)^{-1}
\]
is an isomorphism of $N_{H\times F}(\Delta R)$-sets, where $N_{H\times F}(\Delta R)=P\times(Z\times P)$ acts on
$Z\times P$ by
\[
(u,zv)\cdot x=ux(zv)^{-1}
\]
for every $u,v\in P$, $z\in Z$ and $x\in Z\times P$. Therefore, $k[N_{H\times F}(\Delta R)/\Delta P]e_1\cong k[Z\times P]e_1$ as $k[N_{H\times F}(\Delta R)]$-modules. The latter module is indecomposable, since its restriction to $P\times 1$ is isomorphic to the regular $kP$-module. Moreover, $(k[Z\times P]e_1)[\Delta P]\cong k[Z\times Z(P)]e_1\neq 0$, since the set of $\Delta P$-fixed points of $Z\times P$ is equal to $Z\times Z(P)$ and $e_1$ is $P$-stable. 

On the other hand, as above, since $R=Z(P)$ is characteristic in $P$, the transitivity of the Brauer construction \cite[Proposition~5.8.5]{LinckelmannVol1} gives
\[
(W[\Delta R])[\Delta P]\cong W[\Delta P],
\]
as $k[N_{H\times F}(\Delta P)]$-modules. The latter module is non-zero and indecomposable by Lemma~\ref{lem:candidate_3A6}(b). Therefore, $W[\Delta R]$ has a unique indecomposable direct summand whose Brauer construction at $\Delta P$ is non-zero. By Lemma~\ref{lem:candidate_3A6}(d), $Y$ is isomorphic to a direct
summand of $W[\Delta R]$. We have already shown that $Y$ is indecomposable and that $Y[\Delta P]\neq 0$. Hence, by the preceding paragraph,
\[
Y\cong k[N_{H\times F}(\Delta R)/\Delta P]e_1 \cong k[Z\times P]e_1
\]
as $k[N_{H\times F}(\Delta R)]$-modules. By Remark~\ref{rem:diagonal-normalizers-3A6}, this is an isomorphism of $(kP,k[Z\times P]e_1)$-bimodules. This bimodule induces the required splendid Morita equivalence.
\end{proof}

\vspace{2mm}
\subsection{The splendid Morita equivalences}
Again, it remains to prove that the stable equivalence of Proposition~\ref{prop:gluing-3A6} is a splendid Morita equivalence, which we achieve using~\cite[Theorem 2.1(ii) and~(iii)]{Lin96StableEqs}. We start by collecting information on the simple and trivial source modules of the blocks involved.

\begin{lem}\label{lem:tsmodulesB}\leavevmode
\begin{enumerate}[label={\rm(\alph*)}]
    \item The trivial $kH$-module $k_H$ is the unique trivial source $\bB$-module with vertex $P$ and we have $\chi^{}_{\widehat{k_H}}=\chi^{}_1$\,.
    \item 
    The block $\bB$ contains precisely two pairwise non-isomorphic trivial source modules with vertex $Q$, namely $\Sc(H,Q)$ and a uniserial module $\widetilde{U}:=U(\widetilde{3}_1, \widetilde{3}_2)$ (permuting the labelling of $\tilde{3}_1$ and $\tilde{3}_2$ if necessary) affording the character $\chi^{}_{\widehat{\widetilde{U}}}=\chi^{}_4$\,.
    \item Let $\bb$ denote the Brauer correspondent of $\bB$ in $N_H(Q)\cong \fS_4$ and let $f_{H,Q}$ denote the Green correspondence with respect to $(H,Q,N_H(Q))$ from $H$ to $N_H(Q)$.  Then, $f_{H,Q}(\widetilde{U})=2_{\fS_4}$, the unique two-dimensional $\bb$-module. (See Notation~\ref{nota:Section_3A7}\ref{nota:simplesInkH}.)
\end{enumerate} 
\end{lem}

\noindent Before we proceed to the proofs of Lemma~\ref{lem:tsmodulesB} and the next lemma, 
we note that since $H\leq G$ and $F\leq G$, we may consider again the simple $kG$-module $15$ (respectively the simple $kG$-module~$6$) defined in Notation~\ref{nota:Section_3A7}\ref{nota:3A72-decompMatrixA}, which is a trivial source module with vertex $P$ (respectively $Q$) by Lemma~\ref{lem:tsmInA}.  These allow us to calculate the relevant trivial source modules of $kH$ and $kF$.

\begin{proof}
\begin{enumerate}[label={\rm(\alph*)}] 
    \item It is immediate that $k_H$ is a trivial source module with vertex~$P$ and belongs to $\bB$. Moreover, the number of isomorphism classes of trivial source $kH$-modules with vertex $P$
    is given by the number of $2'$-conjugacy classes of  $N_H(P)/P=P/P\cong \{1\}$, i.e. one. (See e.g. \cite[Theorem~4.6(c)]{LasSurvey}.) The claim follows. 
    \item Similarly, the number of 
     isomorphism classes of trivial source $kH$-modules with vertex $Q$ equals the number of $2'$-conjugacy classes of   $N_H(Q)/Q\cong \fS_4/Q\cong \fS_3$, which is two. Since $kH$ has two blocks, namely $\bB=\bB_0(kH)$ and a block of defect zero, any trivial source $kH$-module with vertex $Q$ belongs to $\bB$.  
     Certainly, $\Sc(H,Q)$ is such a module.
     \par
    Next, by~\cite[II, Lemma 12.6(ii)]{LandrockBook},  $\Res^{G}_{H}(\chi^{}_{\widehat{6}})$ takes non-negative integer values at all $2$-elements. Therefore, it follows from Notation~\ref{nota:Section_3A6}\ref{nota:3A6-decompMatrixB} that the only possibility is 
    \[
        \Res^{G}_{H}(\chi^{}_{\widehat{6}})=\chi^{}_{4}\,.
    \]
    Thus, $\Res^{G}_{H}(6)$ is a trivial source module with vertex $Q$. 
    We then read from the decomposition matrix of $\bB$ that $\Res^{G}_{H}(6)$ is uniserial of length $2$ with composition factors $\widetilde{3}_1$ and $\widetilde{3}_2$. The claim follows.

    \item We have $N_H(Q)\cong \fS_4$ and $N_G(Q)\cong C_3\times \fS_4$. Moreover, the block $\bb$ is splendidly Morita equivalent to the block $\bc$ of $N_G(Q)$ (defined in Notation~\ref{nota:Section_3A7}\ref{nota:3A7_def_c_c3}), see Lemma~\ref{lem:SplendidlyEquivBrauerCorrespondents}(a).
    Thus, if $f$ denotes the Green correspondence with respect to $(G,Q,N_G(Q))$ from $G$ to $N_G(Q)$, then Lemma~\ref{lem:tsmInA}(c) yields 
    \[
        f_{H,Q}(\widetilde{U})=\Res^{N_G(Q)}_{N_H(Q)}(f(6))=\Res^{N_G(Q)}_{N_H(Q)}(2)= 2_{\fS_4}\,.
    \]
\end{enumerate} 
\end{proof}

\begin{lem}\label{lem:tsmodulesA}\leavevmode
\begin{enumerate}[label={\rm(\alph*)}]
    \item The simple $kF$-module $9_1$ is the unique trivial source $\bC$-module with vertex $P$ and we have $\chi^{}_{\widehat{9_1}}=\chi^{}_{17}$\,.
    \item The block $\bC$ contains precisely two pairwise non-isomorphic trivial source modules with vertex $Q$, one of which is the uniserial module $U:=U(3_1,3_3)$ affording the character $\chi^{}_{\widehat{U}}=\chi^{}_{16}$ (permuting the labelling of $3_1$ and $3_3$ if necessary).
    \item Let $\bc$ be the Brauer correspondent of $\bC$ in $N_F(Q)=N_G(Q)\cong C_3\times \fS_4$ and let $f_{F,Q}$ denote the Green correspondence with respect to $(F,Q,N_F(Q))$ from $F$ to $N_F(Q)$.  Then, $f_{F,Q}(U)=2$, i.e. the unique two-dimensional $\bc$-module. (See Lemma~\ref{lem:SplendidlyEquivBrauerCorrespondents}(b).)
\end{enumerate} 
\end{lem}

\begin{proof}
\begin{enumerate}[label={\rm(\alph*)}] 
    \item 
    A GAP calculation \cite{GAP} shows that the $2$-permutation module $W:=\Res^{G}_{F}(15)$ affords the character 
    \[
        \chi^{}_{\widehat{{W}}}=\Res^{G}_{F}(\chi^{}_{\widehat{15}})=\chi^{}_{17}+\chi^{}_{16}.
    \]
Furthermore, \cite[I, 
Theorem 17.3]{LandrockBook} implies that $W$ has $kF$-submodules $U$ and $V$ affording the characters $\chi^{}_{16}$ and $\chi^{}_{17}$, respectively.  We then read from the decomposition matrix of $\bC$ that $V\cong 9_1$ and $U$ is indecomposable with composition factors $3_1$ and $3_3$. It follows that $U\cap V=0$ and so $W=U\oplus V$.  
Hence, $9_1$ is a trivial source module, and $P$ is a vertex of $9_1$ as it is odd-dimensional.
\par
    Recalling that 
    $Q\leq P\leq N_F(P)=N_{G}(P)\leq N_F(Q)=N_{G}(Q)\leq F\leq G$\,, we have that 
     the number of isomorphism classes of trivial source $\bC$-modules with vertex $P$ equals the number of isomorphism classes of trivial source modules with vertex $P$ lying in the Brauer correspondent $\bc$ of $\bC$ (which is also the Brauer correspondent of $\bA$) at the level of $N_F(Q)$, and we have already proved in~Lemma~\ref{lem:SplendidlyEquivBrauerCorrespondents}(c) that this number is one. This proves our claim.
    \item By Lemma~\ref{lem:SplendidlyEquivBrauerCorrespondents}(a) and (c), we also know that the number of isomorphism classes of trivial source $\bC$-modules with vertex $Q$ equals the number of isomorphism classes of trivial source $k\fS_4$-modules with vertex $Q$. This number is  two, as $N_{\fS_4}(Q)/Q\cong \fS_3$, which has two conjugacy classes of $2^\prime$-elements.
    \par
   
    Let $f_{G,F}$ denote the Green correspondence with respect to $(G,P,F)$ from $G$ to $F$. 
    The proof of (a) shows that $f_{G,F}(15)=9_1$, thus, as the Green correspondence and the Brauer correspondence commute, $f_{G,F}(6)$  belongs to $\bC$ as well. It is a trivial source module with vertex $Q$ and, by definition, it is a direct summand of the  $2$-permutation module $W':=\Res^{G}_{F}(6)$\,. Moreover, $\chi^{}_{\widehat{f_{G,F}(6)}}$ is a sum of constituents of $\chi^{}_{\widehat{W'}}=\Res^{G}_{F}(\chi^{}_{\widehat{6}})$
    and by~\cite[II, Lemma 12.6(ii)]{LandrockBook} $\chi^{}_{\widehat{f_{G,F}(6)}}$ takes non-negative integer values at all $2$-elements. Therefore, it follows from Notation~\ref{nota:Section_3A6}\ref{nota:3A6-decompMatrixA} that the only possibility is 
    \[
        \Res^{G}_{F}(\chi^{}_{\widehat{6}})=\chi^{}_{\widehat{f_{G,F}(6)}}=\chi^{}_{16}\,.
    \]
    We then read from the decomposition matrix of $\bC$ that $f_{G,F}(6)$ is uniserial of length $2$ with composition factors $3_1$ and $3_3$. The claim follows. 
    \item  Let $f_{G,Q}$ denote the Green correspondence with respect to $(G,Q,N_{G}(Q))$ from $G$ to $N_G(Q)$. 
    Since $N_F(Q)=N_G(Q)$, we have  
    \[
    f_{F,Q}(U)=f_{F,Q}(f_{G,F}(6))=f_{G,Q}(6)=2
    \]
    by the proof of (b) and Lemma~\ref{lem:tsmInA}(c).
\end{enumerate} 
\end{proof}

\noindent We can now prove that
the functor $M\otimes_{\bC}-$ maps 
the simple $\bC$-modules to the simple \medskip $\bB$-modules.

\begin{lem}\label{lem:simple-tensor-M-3A6}%
Let $M$ be the $(\bB,\bC)$-bimodule defined in Notation~\ref{nota:BimodulesSection_3A6}{\rm\ref{nota:DefBimodule-3A6}}. 
Then, permuting the labelling of $\tilde{3}_1$ and $\tilde{3}_2$ if necessary, the following assertions hold:
\begin{enumerate}[label={\rm(\alph*)}]
    \item $M\otimes_{\bC}9_1\cong k_H$ as $\bB$-modules and  $M^\ast\!\otimes_{\bB}k_H\cong 9_1$ as $\bC$-modules;
    \item $M\otimes_{\bC}U(3_1,3_3)\cong U(\tilde{3}_1,\tilde{3}_2)$ as $\bB$-modules; and 
    \item $M\otimes_{\bC}3_1\cong \tilde{3}_1$ and $M\otimes_{\bC}3_3\cong \tilde{3}_2$  as $\bB$-modules.  
\end{enumerate}
\end{lem}

\begin{proof}
First, by \cite[Theorem 2.1(ii)]{Lin96StableEqs}, the $\bB$-modules $M\otimes_{\bC}9_1$ , $M\otimes_{\bC}3_1$ and $M\otimes_{\bC}3_3$ are all indecomposable and non-projective, since $9_1$, $3_1$ and $3_3$ are non-projective simple $\bC$-modules.
\begin{enumerate}[label={\rm(\alph*)}]
    \item 
     Since $k_H$ is the unique trivial source $\bB$-module with vertex $P$ by Lemma~\ref{lem:tsmodulesB}(a) and $9_1$ is the unique trivial source $\bC$-module with vertex $P$ by Lemma~\ref{lem:tsmodulesA}(a),  it follows from Lemma~\ref {lem:stablesplendid_and_pperm_modules} that 
    $M\otimes_{\bC}9_1\cong k_H$. 
    An analogous argument with $M^\ast$ yields the second claim. 
    \item As above, write $U:=U(3_1,3_3)$ and $\widetilde{U}:=U(\widetilde{3}_1, \widetilde{3}_2)$\,.
    First, by Lemma~\ref{lem:stablesplendid_and_pperm_modules}(b), we have that  $$M\otimes_{\bC}U\cong Y\oplus \mathcal{P}$$
    where $Y$ is a trivial source $\bB$-module with vertex $Q$ and $\mathcal{P}$ is a projective $\bB$-module. 
    It remains to prove that  $Y\cong\widetilde{U}$ and \medskip $\mathcal{P}\cong0$. We proceed in three steps.\\
    (i) First, we claim that $Y=\widetilde{U}$.\\ 
    Recall that $N_H(Q)\cong \fS_4<N_F(Q)\cong C_3\times \fS_4$ and $\Res^{N_F(Q)}_{N_H(Q)}$ induces a splendid Morita equivalence between the Brauer correspondent $\bb$ of $\bB$ in $N_H(Q)$ and the Brauer correspondent $\bc$ of $\bC$ in $N_F(Q)$ (see Lemma~\ref{lem:SplendidlyEquivBrauerCorrespondents}(a)).  Also notice that $X=e_{\bB}\widetilde{X}e_{\bC}\cong e_{\bB}kH\otimes_{kP}kFe_{\bC}=\Res^{H\times H}_{H\times P}(\bB)\otimes_{kP}\Res^{F\times F}_{P\times F}(\bC)$. Thus, we have
    \[
     M\otimes_{\bC}U
     \,\Big|\, \Ind_{N_H(Q)}^H\circ\Res^{N_F(Q)}_{N_H(Q)}\circ\Res^F_{N_F(Q)}(U)
    \]
    and we claim that $\widetilde{U}$ is the unique direct summand with vertex $Q$ of the right-hand side module.
    Considering first  the Green correspondence  with respect to $(F,Q,N_F(Q))$, Lemma~\ref{lem:tsmodulesA}(c) yields 
    \[
       \Res^F_{N_F(Q)}(U) = f_{F,Q}(U)\oplus V = 2\oplus V
    \]
    for some $2$-permutation $kN_F(Q)$-module $V$. Hence, $\Res^F_{N_H(Q)}(U)\cong 2_{\fS_4}\oplus V'$ for some $2$-permutation $kN_H(Q)$-module $V'$ of dimension $4$. 
    Using again~\cite[II, Lemma 12.6(ii)]{LandrockBook} and  the character table of $\fS_4$, we obtain that $V'$ is indecomposable and $\chi^{}_{\widehat{V'}}$ takes value zero at the elements of order 4 and only one of the two conjugacy classes of elements of order $2$, implying that $V'$ has a cyclic vertex of order $2$.
    Next, we consider the Green correspondence with respect to $(H,Q,N_H(Q))$. By Lemma~\ref{lem:tsmodulesB}(c), we have $f_{H,Q}(\widetilde{U})=2_{\fS_4}$, so $\widetilde{U}$ is the unique indecomposable summand of $\Ind_{N_H(Q)}^H(2_{\fS_4})$ with vertex $Q$, and $\Ind_{N_H(Q)}^H(V')$ does not have an indecomposable summand with vertex $Q$. It follows that \medskip $Y=\widetilde{U}$. \\
    (ii) We claim that $P(k_H)\nmid\mathcal{P}$.
    It is enough to prove that $\Hom_{\bB}(M\otimes_{\bC}U, k_H)\cong 0$.
    First, as $k_H$ is simple, by \cite[Corollary 4.13.4]{LinckelmannVol1}, we have  
    $\Hom_{\bB}(M\otimes_{\bC}U, k_H)\cong\,  \underline{\Hom}_{\bB}(M\otimes_{\bC}U, k_H)$. Therefore, it follows from (a) that 
     \begin{align*}
     \Hom_{\bB}(M\otimes_{\bC}U, k_H)
     \cong&\, \underline{\Hom}_{\bB}(M\otimes_{\bC}U, k_H)\\
      \cong&\,  \underline{\Hom}_{\bC}(U, M^\ast\!\otimes_{\bB}k_H)\\
      \cong&\,  \underline{\Hom}_{\bC}(U(3_1,3_3), 9_1)\cong 0\,.
     \end{align*}
    (iii) We claim that $\mathcal{P}\cong 0$.\\
    First recall from Notation~\ref{nota:Section_3A6}\ref{nota:3A6-decompMatrixB} that $\Irr_k(\bB)=\{k_H, \tilde{3}_1,\tilde{3}_2\}$. Hence, it follows from~(ii) that $\mathcal{P}$ is of the form 
    \[
        \mathcal{P}=n_1 P(\tilde{3}_1)\oplus n_2P(\tilde{3}_2)
    \]
   where $n_1,n_2$ are non-negative integers. 
    By Lemma~\ref{lem:tsmodulesB}(b) we have $\chi^{}_{\widehat{\widetilde{U}}}=\chi^{}_4$ and a GAP calculation \cite{GAP} yields
    \[
        \Ind_{N_H(Q)}^H\circ\Res^{N_F(Q)}_{N_H(Q)}\circ\Res^F_{N_F(Q)}(\chi^{}_{\widehat{U}})=\chi^{}_1+3\chi^{}_{4}+\chi^{}_{5}+2\eta
    \]
    where $\eta$ is an irreducible character belonging to the unique block of defect zero of $H$. (See~\cite[L3(2)mod2]{ModularAtlas}.) Thus, we have
    \[
    \chi^{}_{\widehat{\mathcal{P}}}\mid \chi^{}_{1}+2\chi^{}_{4}+\chi^{}_5\,.
    \]
    However, we read from the decomposition matrix of $\bB$ in Notation~\ref{nota:Section_3A6}\ref{nota:3A6-decompMatrixB} that 
    \[
            \chi^{}_{\widehat{P(\tilde{3}_1)}}=\chi^{}_{2}+\chi^{}_{4}+\chi^{}_{5}\qquad\text{ and }\qquad \chi^{}_{\widehat{P(\tilde{3}_2)}}=\chi^{}_{3}+\chi^{}_{4}+\chi^{}_{5}\,,
    \]
    proving that $\mathcal{P}\cong 0$. \\
    \item  
    The bimodule $M$ is left and right projective, so the functor $M\otimes_{\bC}-$ is exact. It follows that the non-split exact sequence 
    \[
        \begin{tikzcd}
            0 \arrow[r]
             & 3_3 \arrow[r, hook, "\iota"]
             & U(3_1,3_3) \arrow[r, two heads, "\pi"]
             & 3_1 \arrow[r]
             & 0
        \end{tikzcd}
    \]
    of $\bC$-modules 
    (where $\iota$ is the inclusion of $3_3$ in the socle of $U(3_1,3_3)$ and $\pi$ the canonical projection onto its head)
    is mapped to the non-split exact sequence 
    \[
        \begin{tikzcd}
            0 \arrow[r]
             & M\otimes_{\bC}3_3 \arrow[r, hook, "M\otimes_{\bC}\iota"]
             & M\otimes_{\bC}U(3_1,3_3) \arrow[r, two heads, "M\otimes_{\bC}\pi"]
             & M\otimes_{\bC}3_1 \arrow[r]
             & 0
        \end{tikzcd}
    \]
    of $\bB$-modules,  
    whose middle term is $M\otimes_{\bC}U(3_1,3_3)\cong U(\tilde{3}_1,\tilde{3}_2)$ by (b).
    Since the latter module has length $2$, the 
    additivity of composition lengths implies that both $M\otimes_{\bC}3_3$ and $M\otimes_{\bC}3_1$ have length $1$, so they are simple and the claim follows. 
\end{enumerate}
\end{proof}

\begin{thm}\label{thm:FaithfulBlocks3A6}%
    The two faithful blocks $\bC,\bC^\ast\in\Bl^{\mathrm{max}}(k[3.\fA_6])$ are splendidly Morita equivalent to the principal block $\bB=\bB_0(k\PSL_2(7))$ of the projective special linear group $\PSL_2(7)$.
\end{thm}

\begin{proof}
By Proposition~\ref{prop:gluing-3A6}, the $(\bB,\bC)$-bimodule $M$ defined in Notation~\ref{nota:BimodulesSection_3A6}{\rm\ref{nota:DefBimodule-3A6}} realises a stable equivalence of Morita type between $\bB$ and $\bC$. Furthermore, by Lemma~\ref{lem:simple-tensor-M-3A6}
for any simple $\bC$-module $S$, the $\bB$-module $M\otimes_{\bC}S$ is simple.  Therefore, $M$ realises a Morita equivalence between $\bB$ and $\bC$ by \cite[Theorem~2.1(iii)]{Lin96StableEqs}.
    Since  $M$ is a  $2$-permutation bimodule by definition, this equivalence is splendid, as required. Replacing $\bC$ with $\bC^\ast$, a completely analogous argument proves that there is a splendid Morita equivalence between $\bC^{\ast}$ and $\bB$. 
\end{proof}



\section{Proof of the main results}\label{sec:proof_Main_Thm}\label{sec:proofs}%

Finally, we prove our main result.

\begin{proof}[Proof of Theorem~\ref{thm:mainthm_intro}]
Let $\bB\in\Bl^{\mathrm{max}}(kG)$. By Lemma~\ref{lem:Fong2ndredDihedral}, $\bB$ is splendidly Morita equivalent to a block $\widehat{\bB}$ of a central extension $\widehat{L}$  of a group \(L\), with dihedral Sylow \(2\)-subgroups of order \(2^n\) and \(O_{2'}(L)=1\), by a cyclic group of order $|L|_{2'}$. Clearly, $\widehat{\bB}\in\Bl^{\mathrm{max}}(k\widehat{L})$, and by Subsection~\ref{ssec:dihedralSylow2} the group $L$ is of type {\bf(D1)}, {\bf(D2)}, {\bf(D3)} or {\bf(D4)}. 
\par
Now, first assume that $L\notin\{\mathfrak{A}_7, \PSL_2(9)\}$. Then, $H^2(L,k^\times)$ is trivial by Lemma~\ref{lem:SchurMult_2prime}. Therefore, we can take $\widehat{L}=L$ by Theorem~\ref{thm:Fong2ndReduction}(c), and the claim follows from Proposition~\ref{prop:ClassificBlocksDihedralDefectGroupsMaxDefectO_2prime_is_1} in this case. 
\par 
Next, assume that  $L=\PSL_2(9)\cong\fA_6$. 
Let $Z$ be the cyclic
group of order $|L|_{2'}=45$. By Lemma~\ref{lem:SchurMult_2prime}, we have
\[
{\mathrm H}^2(L,k^\times)\cong C_3.
\]
It follows that, besides the trivial extension $Z\times L$, 
up to isomorphism of groups, there is a unique non-split central
extension $H$ of $L$ by $Z$. This extension is obtained by starting from the Schur
cover
\[
1\longrightarrow C_3\longrightarrow 3.L\longrightarrow L\longrightarrow 1
\]
and pushing out along the embedding $C_3\hookrightarrow Z$. Thus, if we write
$\widetilde{L}:=3.L$, then
\[
H\cong (\widetilde{L}\times Z)/\Delta(C_3),
\]
that is, $H$ is the central product of $\widetilde{L}$ and $Z$, obtained by
identifying the central subgroup of order $3$ in $\widetilde{L}$ with the unique
subgroup of order $3$ in $Z$.
\par
If $\widehat{L}=Z\times L$, then any block $\widehat{\bB}\in\Bl^{\mathrm{max}}(k\widehat{L})$ is splendidly Morita equivalent to $\bB_0(kL)$ by Lemma~\ref{lem:SplendidMoritaEquivalence_direct_and_central_product}, because the principal block of $L$ is the unique block of maximal defect of $L$ by Proposition~\ref{prop:ClassificBlocksDihedralDefectGroupsMaxDefectO_2prime_is_1}.
Now, if $\widehat{L}=H$, 
by \cite[Proposition~1.5]{Duvel04} (see Subsection~\ref{ssec:blocks_of_central_producrs}), every block $\widehat{\bB}\in\Bl^{\mathrm{max}}(kH)$ is splendidly Morita equivalent to  $\widetilde{\bB}\otimes_k kZe_\lambda$, where $\widetilde{\bB}$ is a
maximal-defect block of $k\widetilde{L}$ and $e_\lambda$ is a block idempotent of
the split semisimple algebra $kZ$ corresponding to an irreducible character
$\lambda$ of $Z$.  Moreover, each block
$\widetilde{\bB}\otimes_k kZe_\lambda$ is splendidly Morita equivalent to
$\widetilde{\bB}$ by Lemma~\ref{lem:SplendidMoritaEquivalence_direct_and_central_product}.
\par
For the case $L=\fA_7$, we have $|L|_{2'}=315$ and  an analogous  argument works replacing the cyclic group $C_{45}$ with $C_{315}$.
\par 
Therefore, it remains to consider the blocks $\widetilde{\bB}\in \Bl^{\textrm{max}}(k\widetilde{L})$ for $\widetilde{L}\in\{3.\fA_6,3.\fA_7\}$. First, if $\widetilde{\bB}$ is the principal block, then there is nothing to do. Indeed, since the centre has order $3$, inflation from $\widetilde{L}/Z(\widetilde{L})\cong L$ to  $\widetilde{L}$ induces a splendid Morita equivalence between 
$\bB_0(kL)$ and $\bB_0(k\widetilde{L})$.
Finally, if $\widetilde{B}$ is a faithful block of $k\widetilde{L}$, then $\widetilde{\bB}$ is splendidly Morita equivalent to $\bB_0(k\PSL_2(7))$ by Theorem~\ref{thm:FaithfulBlocks3A6} provided $\widetilde{L}=3.\fA_6$, whereas $\widetilde{\bB}$ is splendidly Morita equivalent to $\bB_0(k\fS_4)$ by Theorem~\ ref{thm:FaithfulBlocks3A7} provided $\widetilde{L}=3.\fA_7$. Since $\fS_4\cong\PGL_2(3)$, the assertion is proved.
This gives us all the blocks listed in cases {\rm(a)--(f)}.
\par
In cases {\rm(c)--(f)}, any choice of $q$ satisfying the indicated condition yields a representative of the corresponding splendid Morita equivalence class by Proposition~\ref{prop:ClassificBlocksDihedralDefectGroupsMaxDefectO_2prime_is_1}. 
\end{proof}

%% file: appendix.tex

\appendix
\section{Proof of the vertex and source statement in Fong's second reduction}\label{app:A}%

In this Appendix, we give a proof of Part (b) of Theorem~\ref{thm:Fong2ndReduction} establishing that Fong's second reduction can be understood as a Morita equivalence with endo-permutation source. This result, understood in terms of source-algebras, is due to Puig and we refer to \cite[Theorem~9.11.9] {LinckelmannVol2}, which states the connection between the source-algebras of the blocks involved and the source of the bimodule giving the Morita equivalence.  Below we give a precise proof exhibiting a source and a vertex of this bimodule purely in module-theoretic terms, which in particular does not require knowledge of the source-algebras.

\begin{proof}[Proof of Theorem~\ref{thm:Fong2ndReduction}, Part {\rm(b)} and
Remark~\ref{rem:ProofFong2ndRed}]{\leavevmode}

\smallskip
\noindent
{\rm(i)} \emph{The module $V$ and its restriction to $D$.}
Since $N\cap D=1$, the quotient group $ND/N\cong D$ is a $p$-group. Moreover, since $\bc$ is $G$-stable, the inertia group of the simple $\bc$-module $V$ in $ND$ is $ND$. Therefore, by
\cite[Chapter~III, Corollary~3.16]{FeitBook}, $V$ has a unique extension, up to isomorphism, to an irreducible $k[ND]$-module. We denote this extension again by $V$ and set
\( W:=\Res_D^{ND}(V)\). 
Consider $S=kNe_{\bc}$ as a $kD$-module via conjugation. Then
\[
S=kNe_{\bc}\cong\End_k(V)\cong\End_k(W)
\]
as $kD$-modules. Since $S$ is a direct summand of the permutation $kD$-module $kN$, it follows that $\End_k(W)$ is a permutation $kD$-module. Hence $W$ is an endo-permutation $kD$-module.

Moreover, $\Br_D(e_{\bc})\neq0$, and hence $\End_k(V)\cong kNe_{\bc}$ has a trivial direct summand as a $kD$-module. It follows from \cite[Lemma~2.2]{ThevenazSurvey} that $W$ is capped. We set
\[
U:=\operatorname{cap}(W).
\]
Thus $U$ is an indecomposable $kD$-module with vertex $D$.

We next show that the $kND$-module $V$ has vertex $D$. Let $Q$ be a vertex of $V$. Since $\bc$ is a block of defect zero, the module $V$ belongs to a block of $kND$ having a defect group $P$ such that $P\cap N=1$. Therefore, $|Q|\leq |P|\leq |D|$. On the other hand, $U$ is a direct summand of $\Res_D^{ND}(V)$ and has vertex $D$. Hence $D$ is contained, up to conjugation, in a vertex of $V$. It
follows that $|D|\leq |Q|$. Therefore $|Q|=|D|$, and $D$ is a vertex of $V$.

\smallskip
\noindent
{\rm(ii)} \emph{Construction of the Morita equivalence.}
By \cite[Section~3]{SalminenSOURCES08}, there exists an isomorphism
\[
\Phi:S\otimes_k kHe\xrightarrow{\sim}kGe_{\bc}
\]
of interior $D$-algebras, where $e$ is an idempotent of $kZ$. Here, we identify the copy of $D$ in $H$ arising from this construction with $D$.

Since $S$ is simple, there is a bijection between the block idempotents of $kH$ appearing in $e$ and the blocks of $kG$ appearing in $e_{\bc}$; see
\cite[Section~4]{SalminenSOURCES08}. Let
$\widehat{\bB}$ be the block of $kH$ corresponding to
$\bB$. By the proof of
\cite[Proposition~4.1]{SalminenSOURCES08}, which does not use the nilpotency assumption, the chosen copy of $D$ in $H$ is a defect group of $\widehat{\bB}$.

The isomorphism $\Phi$ restricts to an isomorphism $\Phi_{\bB}:S\otimes_k\widehat{\bB}\xrightarrow{\sim}\bB$ of interior $D$-algebras. Since $S\cong\End_k(V)$, the $(S\otimes_k\widehat{\bB},\widehat{\bB})$-bimodule $V\otimes_k\widehat{\bB}$ induces a Morita equivalence between
$S\otimes_k\widehat{\bB}$ and $\widehat{\bB}$. Its actions are given by
\[
(s\otimes a)\cdot(v\otimes x)\cdot a'=s(v)\otimes axa'
\]
for all $s\in S$, $a,a',x\in\widehat{\bB}$ and $v\in V$; see \cite[Exercise~9.5(a)]{ThevenazBook}. Transporting the left action along $\Phi_{\bB}^{-1}$, the bimodule $M:=V\otimes_k\widehat{\bB}$ becomes a $(\bB,\widehat{\bB})$-bimodule and induces a Morita equivalence between $\bB$ and $\widehat{\bB}$.

\smallskip
\noindent
{\rm(iii)} \emph{Relative projectivity of $M$.}
Set $\overline G:=G/N$
and let $\pi:H\to\overline G$
be the canonical surjection. Consider the $(kGe_{\bc},kHe)$-bimodule $\mathcal M:=V\otimes_k kHe$, where the left action is transported along $\Phi:S\otimes_k kHe\xrightarrow{\sim}kGe_{\bc}$. 
Recall that $Z$ is identified with a finite subgroup of $k^\times$ and that
\[
e=\frac{1}{|Z|}
\sum_{\zeta\in Z}\zeta^{-1}(1,\zeta).
\]
For every $\zeta_0\in Z$, we have $(1,\zeta_0)e=\zeta_0e$, where, on the right-hand side, $\zeta_0$ is regarded as an element of $k^\times$. For every $\overline x\in\overline G$, choose an element $\widetilde x\in H$ with $\pi(\widetilde x)=\overline x$. Since $Z$ is central in $H$, we have
\[
(1,\zeta)\widetilde xe=\widetilde x(1,\zeta)e=\zeta\widetilde xe \qquad(\zeta\in Z).
\]
It follows that $kZ\widetilde xe=k\widetilde xe$. Therefore, as a $k$-vector space,
\[
kHe=\bigoplus_{\overline x\in\overline G} k\widetilde xe,
\]
and hence
\[
\mathcal M=\bigoplus_{\overline x\in\overline G}
\mathcal M_{\overline x}\;\;\text{ with }\quad\mathcal M_{\overline x}:=V\otimes_k k\widetilde xe.
\]
We now use the explicit construction of the same
isomorphism $\Phi$. By \cite[Section~3]{SalminenSOURCES08}, for every $g\in G$
there exist $s_g\in S^\times$ and $\alpha_g\in k^\times$ such that
\[
\Phi^{-1}(ge_{\bc})=s_g\otimes \alpha_g\widetilde{gN}e,
\tag{1}
\]
where $\widetilde{gN}\in H$ is a lift of
$gN\in\overline G$. 
Let $\overline x\in\overline G$, $v\otimes y\in\mathcal M_{\overline x}$ and $y\in k\widetilde xe$. Using~{\rm(1)}, we obtain
\begin{align*}
(g,h)\cdot(v\otimes y)&=g\cdot(v\otimes y)\cdot h^{-1}=\alpha_gs_g(v)\otimes \widetilde{gN}yh^{-1}.
\end{align*}
The image of
$\widetilde{gN}\widetilde xh^{-1}$ under $\pi$ is
$(gN)\overline x\pi(h)^{-1}$. Therefore $(g,h)\mathcal M_{\overline x} \subseteq \mathcal M_{(gN)\overline x\pi(h)^{-1}}$.
Since the map induced by $(g,h)$ is invertible and both spaces have dimension $\dim_k(V)$, we obtain
\[
(g,h)\mathcal M_{\overline x}=\mathcal M_{(gN)\overline x\pi(h)^{-1}}.
\tag{2}
\]
In particular, $(g,h)\mathcal M_{\overline1}=\mathcal M_{\overline1}$ if and only if $gN=\pi(h)$. Thus the stabilizer of $\mathcal M_{\overline1}=V\otimes_k ke$ in $G\times H$ is
\[
L:=\left\{(g,h)\in G\times H\ \middle|\ gN=\pi(h) \right\}.
\]
Moreover, Equation~{\rm(2)} shows that $G\times H$ acts transitively on $\left\{\mathcal M_{\overline x}\ \middle|\ \overline x\in\overline G\right\}$. It follows that the natural homomorphism
\[
\Theta:\Ind_L^{G\times H}(\mathcal M_{\overline1})\longrightarrow\mathcal M,\qquad(g,h)\otimes m\longmapsto(g,h)\cdot m,
\]
is surjective. Since both sides have the same dimension, $\Theta$ is an isomorphism, so $\mathcal M$ is induced from $\mathcal M_{\overline1}$. Since $\widehat{\bB}$ is a block summand of $kHe$, the bimodule $M$ is a direct summand of $\mathcal M$ and is therefore
relatively $L$-projective.

\smallskip
\noindent
{\rm(iv)} \emph{Vertex and source of $M$.}
Since $M$ induces a Morita equivalence between the block algebras $\bB$ and $\widehat{\bB}$, it is indecomposable as a $k[G\times H]$-module. Let $R$ be a vertex of $M$. Since $M$ is relatively $L$-projective, $R$ is contained in a conjugate of $L$. On the other hand, $M$ belongs to the block $\bB\otimes_k\widehat{\bB}^*$ of $k[G\times H]$. Hence, after replacing $R$ by a conjugate, we may assume that $R\leq (D\times D)\cap{}^yL$ for some $y\in G\times H$.

Since $D\cap N=1$ and $D\cap Z=1$, it follows from the definition of $L$ that $(D\times D)\cap{}^yL$ is a twisted diagonal
subgroup of $D\times D$.
In particular, $|R|\leq |D|$.

Viewed as a $k[H\times H]$-module, the regular bimodule $\widehat{\bB}$ has vertex $\Delta D$ and trivial source. Therefore $k\mid
\Res_{\Delta D}^{H\times H}(\widehat{\bB})$. Since $\Phi_{\bB}$ is an isomorphism of interior
$D$-algebras, the action of $\Delta D$ on $M$ is given by
\[
(d,d)\cdot(v\otimes x) =dv\otimes dxd^{-1}.
\]
Therefore,
\[
\Res_{\Delta D}^{G\times H}(M)\cong W\otimes_k
\Res_{\Delta D}^{H\times H}(\widehat{\bB})
\]
as $k[\Delta D]$-modules, where we identify
$D\cong\Delta D$. 
It follows that $W\mid\Res_{\Delta D}^{G\times H}(M)$. Since $U\mid W$, we obtain $U\mid \Res_{\Delta D}^{G\times H}(M)$.
Viewed as a $k[\Delta D]$-module, $U$ has vertex
$\Delta D$. Therefore, there exists $z\in G\times H$ such that $\Delta D\leq{}^zR$.
Hence $|D|\leq |R|$. Together with the preceding inequality, this implies that $|R|=|D|$ and ${}^zR=\Delta D$. Thus, $\Delta D$ is a vertex of $M$.

Finally, the indecomposable $k[\Delta D]$-module $U$ is a direct summand of $\Res_{\Delta D}^{G\times H}(M)$ and has vertex $\Delta D$. Therefore, $U$ is a source of
$M$. In other words, via the identification
$D\cong\Delta D$, the cap of $\Res_D^{ND}(V)$ is a source of the bimodule inducing the Morita equivalence.
\end{proof}